\newcommand{\CLASSINPUTinnersidemargin}{17mm}
\documentclass[journal]{IEEEtran}
\ifCLASSINFOpdf \else \fi
\usepackage{amssymb,latexsym}
\usepackage{amsfonts}
\usepackage{amsthm,graphicx,epsfig}
\usepackage[cmex10]{amsmath}
\usepackage{cite}

\numberwithin{equation}{section}
\renewcommand{\theequation}{\thesection.\arabic{equation}}
\newtheorem{Thm}{Theorem}
\newtheorem{Lem}{Lemma}

\renewcommand{\baselinestretch}{1.1}

\begin{document}
\title{Consistent estimation of non-bandlimited spectral density from uniformly spaced samples}
\author{Radhendushka Srivastava and~Debasis Sengupta~
\thanks{The authors are with the Applied Statistics Unit, Indian Statistical Institute, Kolkata,
700108, India (e-mail: radhe\_r@isical.ac.in;
sdebasis@isical.ac.in)}}

\maketitle

\maketitle

\begin{abstract}
In the matter of selection of sample time points for the estimation
of the power spectral density of a continuous time stationary
stochastic process, irregular sampling schemes such as Poisson
sampling are often preferred over regular (uniform) sampling. A
major reason for this preference is the well-known problem of
inconsistency of estimators based on regular sampling, when the
underlying power spectral density is not bandlimited. It is argued
in this paper that, in consideration of a large sample property like
consistency, it is natural to allow the sampling rate to go to
infinity as the sample size goes to infinity. Through appropriate
asymptotic calculations under this scenario, it is shown that the
smoothed periodogram based on regularly spaced data is a consistent
estimator of the spectral density, even when the latter is not
band-limited. It transpires that, under similar assumptions, the
estimators based on uniformly sampled and Poisson-sampled data have
about the same rate of convergence. Apart from providing this
reassuring message, the paper also gives a guideline for
practitioners regarding appropriate choice of the sampling rate.
Theoretical calculations for large samples and Monte-Carlo
simulations for small samples indicate that the smoothed periodogram
based on uniformly sampled data have less variance and more bias
than its counterpart based on Poisson sampled data.
\end{abstract}


\section{Introduction}

Estimation of the power spectral density of a continuous time wide
sense stationary stochastic process is an old problem. A set of
regularly (uniformly) spaced samples is generally used for this
purpose. When the process is bandlimited, the spectral density of
the original process can be recovered from that of the sampled
process, provided that the sampling is fast enough. In such a case,
estimation of the spectral density from finitely many observations
at an appropriate sampling rate is a well established topic and many
useful nonparametric and parametric methods have been developed
\cite{Kay}. If the underlying process is not bandlimited, the
spectral density of the original process is not identifiable from
regularly spaced samples, because of the problem of wrapping around
of the spectral density caused by the process of sampling -- also
known as aliasing \cite{Shannon}. In such a case, one cannot
estimate the spectral density consistently from regularly spaced
samples at any fixed sampling rate.

This inadequacy of sampling at regular intervals necessitated the
exploration of other strategies for sampling. Shapiro and Silverman
\cite{Silverman} considered alias free sampling schemes in the sense
that two continuous time processes with different power spectra do
not produce the same spectrum of the sampled sequence. They proved
that additive random sampling through a class of renewal processes,
including the homogeneous Poisson process, is alias free. Beutler
\cite{Beutler} formalized the definition of alias free sampling
relative to a family of spectral distributions and gave different
sampling schemes that are alias free relative to different families
of spectral distributions. Masry \cite{Masrynew} provided a modified
definition of alias free sampling that would guarantee the existence
of a consistent estimator of the power spectral density.

Subsequently, the possibility of breaking free from the nuisance of
aliasing through irregular sampling enthused many researchers and
practitioners. Masry \cite{Masrypoisson} proposed a Poisson sampling
based estimator similar to the smoothed periodogram, and proved that
the proposed estimator is consistent for any average sampling rate,
under certain conditions. Several irregular sampling based
methodologies of estimation of the power spectral density of a
non-bandlimited process have also been proposed
\cite{Mitchel}--\cite{Stoica2009} and applied to various fields
including signal and image processing. Some of these methods are
analogous to methods developed for regularly sampled data.

Irregularly spaced data occur naturally in many practical situations
like seismic data \cite{Ferber}, turbulent velocity fluctuation
\cite{Tummers}, Laser Doppler Anemometer (LAD) data \cite{Nobach},
wide-band antenna arrays \cite{Ishimaru}, Computer Aided Tomography,
Spotlight-Mode Synthetic Aperture Radar \cite{Munson} and so on.
There have been attempts to use standard methods for spectrum
estimation for such data, after suitable weighting\cite{Bronez} or
interpolation \cite{Tummers} of the data. The advent of new
methodologies for irregularly sampled data are important for such
problems.

However, in many applications, including internet traffic data
\cite{Roughan}, seismology \cite{Hung}, \cite{Costain}, image
processing \cite{Eldar} and so on, one has control over the sampling
mechanism. In such cases, the selection of the sampling scheme is a
serious issue. Regular sampling is easier to implement than
irregular sampling. The literature on power spectrum estimation
based on regular sampling contains a large collection of methods,
and these have been studied in detail. An estimator based on
regularly sampled data is generally computationally simpler than a
similar estimator based on irregularly sampled data. Moreover, it is
well known that spectral estimators based on irregularly sampled
data have higher variances than those based on regularly sampled
data \cite{Robert}, \cite{Moore}. On the other hand, the possibility
of aliasing and the proven inconsistency of spectral estimators
based on regular sampling are arguments in favour of using irregular
sampling.

These opposing arguments necessitate a thorough comparison of
spectral estimators based on different sampling strategies. A systematic numerical
comparison is not available in the literature. A
simulated example given by Masry \cite{Masryperformance} was meant
only to highlight the problem of aliasing and to demonstrate how it
can be overcome with irregular sampling. Studies by Roughan
\cite{Roughan} in the special case of active measurements for
network performance produced mixed results, which led the author to
conclude that, while spectral estimators based on Poisson sampling
have less efficiency (i.e., high variance), such techniques could be
used to detect periodicities in the system, and to determine which
rate of regular sampling would be inadequate. In the absence of a
comprehensive empirical study, it appears that many researchers shun
regular sampling mainly because of the stigma of inconsistency
attached to spectral estimators based on regularly sampled data
\cite{Wolf}.

Consistency of an estimator concerns its behaviour as the sample
size goes to infinity. However, it does not make practical sense to
let the sample size tend to infinity with fixed sampling rate. If
one gathers more and more resources to increase the sample size, one
can use some of these resources to sample faster. Realizing this,
practitioners fix the intended range of spectrum estimation, and
then sample an appropriately filtered process at a sufficiently high
frequency to avoid aliasing. Sometimes one goes for successively
higher rates of uniform sampling to determine an appropriate rate of
sampling \cite{Eldar}. However, these common sense approaches are
yet to be backed up by appropriate asymptotic calculations. There is
a need to bridge this gap by working out the large sample properties
of estimators {\it when the sampling rate changes suitably as the
sample size goes to infinity}.

Let $\phi(\cdot)$ be the spectral density of a mean square
continuous, wide sense stationary process $\{X(t),\
-\infty<t<\infty\}$. The most commonly used nonparametric estimator
of $\phi(\cdot)$ based on the $n$ uniformly spaced samples at the
sampling rate $\rho$ is $\widehat{\phi}_{\rho}(\cdot)$, where
\begin{equation}\label{unifest}
\widehat{\phi}_{\rho}(\lambda)=\frac{1}{2\pi\rho}\sum_{|v|<n}\widehat{\gamma}_\rho(v)K(b_n
v)e^{-\frac{iv\lambda}{\rho}}%
1_{[-\pi\rho,\pi\rho]}(\lambda),
\end{equation}
where $K(\cdot)$ is a covariance averaging kernel, $1_A(\lambda)$ is
the indicator function that takes the value 1 when $\lambda\in A$
and the value 0 otherwise, $\{b_n\}$ is a sequence of window widths
such that $b_n\rightarrow0$ and $nb_n\rightarrow\infty$ as
${n\rightarrow\infty}$ and $\widehat{\gamma}_\rho(v)$ is an
estimator of the covariance function, defined as
\begin{equation*}
\begin{split}
\widehat{\gamma}_\rho(v)=&\frac{1}{n}\sum_{j=1}^{n-|v|}X\left(\frac{j}{\rho}\right)X\left(\frac{j+|v|}{\rho}\right),\\
&\hskip.6in v=-(n-1),-(n-2),\ldots,n-2,n-1. %
\end{split}
\end{equation*}
This estimator is known to be consistent when the underlying process
is bandlimited.

In the present work, without assuming that the process is
bandlimited, we examine the asymptotic properties of the estimator
given in (\ref{unifest}) by letting $\rho$ go to infinity at an
appropriate rate, as $n$ goes to infinity. In the sequel, we shall
use the notation $\rho_n$ instead of $\rho$, in order to explicitly
indicate the dependence of the sampling rate on the sample size.
Accordingly, we use the following modified notation of the estimator
of (\ref{unifest}):
\begin{equation}\label{unifestn}
\widehat{\phi}_n(\lambda)=%
\frac{1}{2\pi\rho_n}\sum_{|v|<n}\widehat{\gamma}_n(v)K(b_n
v)e^{-\frac{iv\lambda}{\rho_n}} %
1_{[-\pi\rho_n,\pi\rho_n]}(\lambda)
\end{equation}
where
\begin{equation}
\begin{split}
\widehat{\gamma}_n(v)=&\frac{1}{n}\sum_{j=1}^{n-|v|}X\left(\frac{j}{\rho_n}\right)X\left(\frac{j+|v|}{\rho_n}\right),\\
&\hskip.6in v=-(n-1),-(n-2),\ldots,n-2,n-1. %
\end{split}
\label{as}
\end{equation}

In Section~\ref{s2}, we prove the consistency of the estimator
$\widehat{\phi}_n(\cdot)$ under some general conditions. In
Section~\ref{s3}, we calculate the rate of convergence of the bias
and variance of this estimator, and determine the optimal rates at
which $\rho_n$ and $nb_n$ should go to infinity so that the mean
square error has the fastest possible rate of convergence.
Subsequently, we compare the rates of convergence of the bias and
the mean squared error (MSE) of this estimator with those of a
similar estimator proposed by Masry \cite{Masrypoisson}, based on
non-uniform (Poisson) sampling. We present the results of a
simulation study in Section~\ref{s4} and provide some concluding
remarks in Section~\ref{s5}. The proofs are given in the appendix.

\section{Consistency of the estimator}\label{s2}
Consider the mean square continuous, wide sense stationary
stochastic process $\{X(t),~ -\infty<t<\infty\} $ with zero mean,
(auto-)covariance function $C(\cdot)$ and spectral density
$\phi(\cdot)$. In order to prove that the estimator (\ref{unifestn})
is consistent, it is sufficient to show that the bias and the
variance of the estimator tend to zero as the sample size ($n$)
tends to infinity.

We assume the following condition on the covariance function
$C(\cdot)$.

\bigskip\noindent
{\bf Condition 1.} The function $h_0(\cdot)$, defined over the real
line as $h_0(t)=\sup_{s\ge|t|}|C(s)|$ is integrable.

\bigskip\noindent
{\bf Remark 1.} Condition 1 is equivalent to saying that the
covariance function $C(\cdot)$ is bounded over $[0,\infty)$ by a
non-negative, non-increasing and integrable function.

\bigskip
We assume the following conditions on the choice of the kernel
$K(\cdot)$, the kernel window width $b_n$ and the sampling rate
$\rho_n$.

\bigskip\noindent
{\bf Condition 2.} The covariance averaging kernel function
$K(\cdot)$ is continuous, even,
square integrable and bounded by a non-negative, even and integrable
function having a unique maximum at 0.
Further, $K(0)=1$.

\bigskip\noindent
{\bf Condition 3.} The kernel window width is such that
$b_n\rightarrow 0$ and $nb_n\rightarrow\infty$ as $n\rightarrow
\infty$.

\bigskip\noindent
{\bf Condition 4.} The sampling rate is such that
$\rho_n\rightarrow\infty$ and $\rho_n b_n\rightarrow 0$ as
$n\rightarrow \infty$.

\bigskip\noindent
{\bf Remark 2.} Note that the estimator $\widehat{\phi}_n(\cdot)$
can be written as
$$\widehat{\phi}_n(\lambda)=
{1\over \rho_nb_n}\int_{-\infty}^\infty
I_n(\mu)W_n\left({\lambda-\mu\over \rho_nb_n}\right)d\mu,$$ where
\begin{eqnarray*}
I_n(\lambda)&=&%
\frac1{2\pi n\rho_n}\left|\sum_{t=1}^nX\!\left(\frac
t{\rho_n}\right)e^{-\frac{it\lambda}{\rho_n}}\right|^21_{[-\pi\rho_n,\pi\rho_n]}(\lambda),\\
W_n(\lambda)&=&\sum_{j=-\infty}^\infty \int_{-\infty}^\infty
K(t)e^{-\frac{it\lambda}{\rho_n}-2\pi i jt}dt.
\end{eqnarray*} %
Thus, $\widehat{\phi}_n(\cdot)$ is the smoothed version of
$I_n(\cdot)$, the periodogram, where the degree of smoothness is
controlled by the smoothing parameter $\rho_nb_n$ of the frequency
domain window $W_n(\cdot)$. Condition 4 says that this parameter
goes to zero, and the sampling rate goes to infinity, as the sample
size goes to infinity.

\bigskip
\begin{Thm}
Under Conditions 1--4, the bias of the estimator
$\widehat{\phi}_n(\cdot)$ given by (\ref{unifestn}) tends to zero
uniformly over any closed and finite interval.
\end{Thm}

Before examining the variance of the estimator we assume the
following condition on the fourth order moments of the process
$\{X(t),~ -\infty<t<\infty\}$.

\bigskip\noindent
{\bf Condition 5.} The fourth moment $E|X(t)|^{4}$ exists for every
$t$, and the fourth moment function
$E[X(t)X(t+v_1)X(t+v_2)X(t+v_3)]$ is a function only of the lags
$v_1$, $v_2$ and $v_3$. Further, the fourth order cumulant function
$Q(v_1,v_2,v_3)$, defined by
\begin{IEEEeqnarray*}{rCl}
Q(v_1,v_2,v_3)&=&P(v_1,v_2,v_3)-P_{G}(v_1,v_2,v_3),\\
{\mbox{where }}\hskip50pt\mbox{}&&\\
P(v_1,v_2,v_3)&=&E(X(t)X(t+v_1))X(t+v_2)X(t+v_3)\\
\mbox{and }P_{G}(v_1,v_2,v_3)
&=&C(v_1)C(v_2-v_1)+C(v_2)C(v_3-v_1)\\&&+C(v_3)C(v_1-v_2),%
\end{IEEEeqnarray*}
satisfies
$$|Q(v_1,v_2,v_3)|\le \prod_{i=1}^{3}g_{i}(v_i),$$
where $g_{i}(v),\ i=1,2,3,$ are all continuous, even, nonnegative
and integrable functions over the real line, which are
non-increasing over $[0,\infty)$.

\bigskip\noindent
{\bf Remark 3.} Condition 5 is satisfied by a Gaussian process, as
the function $P(\cdot)$ reduces to the function $P_G(\cdot)$.

\bigskip
\begin{Thm}
Under Conditions 1--5, the variance of the estimator
$\widehat{\phi}_n(\cdot)$ given by (\ref{unifestn}) converges as
follows:
\begin{equation*}
\lim_{n\rightarrow\infty}nb_nVar[\widehat{\phi}_n(\lambda)]
=
(1+\delta_{0,\lambda})[\phi(\lambda)]^{2}\int_{-\infty}^{\infty}K^{2}(x)dx,
\end{equation*}
where $\delta_{0,\lambda}$ is 1 if $\lambda=0$ and is 0 otherwise.
The convergence is uniform over any closed and finite interval that
does not include the frequency 0. In particular, the variance
converges to 0.
\end{Thm}

\bigskip
It follows from Theorems 1 and 2 that, under Conditions 1--5, the
estimator $\widehat{\phi}_n(\lambda)$ is consistent, and is
uniformly consistent over any closed and finite frequency interval
that does not include the point 0.

\section{Rates of Convergence}\label{s3}
The rate of convergence of the variance of
$\widehat{\phi}_n(\lambda)$ follows from Theorem~2. We assume a few
further conditions in order to arrive at a rate of convergence for
its bias. These include additional conditions on the shapes of the
covariance function and the kernel function.

\bigskip\noindent
{\bf Condition 1A.} The function $h_q(\cdot)$, defined over the real
line as $h_{q}(t)=\sup_{s\ge |t|}|s|^{q}|C(s)|$
 is integrable, for some positive number $q$ greater than 1.

\bigskip\noindent
{\bf Condition 1B.} The spectral density is such that, for some
$p>1$, $\phi(\lambda)$ is $O(|\lambda|^{-p})$, i.e.,
$\displaystyle\lim_{\lambda\rightarrow\infty}|\lambda|^{p}\phi(\lambda)=A$
for some positive $A$.

\bigskip

For any kernel $K(\cdot)$, let us define
$$k_r=\lim_{x\rightarrow 0}\frac{1-K(x)}{|x|^r}$$
for each positive number $r$ such that the limit exists. The
characteristic exponent of the kernel is defined as the largest
number $r$, such that the limit exists and is non-zero
\cite{Parzen}. In other words, the characteristic exponent is the
number $r$ such that $1\!-\!K(1/y)$ is $O(y^{-r})$.

\bigskip\noindent
{\bf Condition 2A.} The characteristic exponent of the kernel $K(\cdot)$
is equal to the number $q$, for which Condition 1A is assumed to
hold.

\bigskip\noindent
{\bf Remark 4.} Condition 1A implies Condition 1 (see Remark 1), and
also that $\phi(\cdot)$ is $[q]$ times differentiable, where $[q]$
is the integer part of $q$. Thus, the number $q$ indicates the
degree of smoothness of the spectral density. If Condition 1A holds
for a particular value of $q$, then it would also hold for smaller
values.

\bigskip\noindent
{\bf Remark 5.} The number $p$ indicates the rate of decay of the
spectral density. The following are two interesting situations,
where Condition~1B holds.
\begin{enumerate}
\item The spectral density $\phi(\cdot)$ is a
rational function, i.e.,
$\phi(\lambda)=\frac{P(\lambda)}{Q(\lambda)}$, where $P(\cdot)$ and
$Q(\cdot)$ are polynomials such that the degree of $Q(\cdot)$ is
more than degree of $P(\cdot)$ by at least $p$. Note that continuous
time ARMA processes possess rational power spectral density.
\item The function $C(\cdot)$ has the following smoothness property:
$C(\cdot)$ is $p$ times differentiable and the $p^{\rm th}$
derivative of $C(\cdot)$ is in $L^1$.
\end{enumerate}

\bigskip\noindent
{\bf Remark 6.} The number $p$ can be increased indefinitely by
continuous time low pass filtering with a cut off frequency larger
than the maximum frequency of interest. There are well-known filters
such as the Butterworth filter, which have polynomial rate of decay
of the transfer function with specified degree of the polynomial,
that can be used for this purpose.

\bigskip
\begin{Thm}
Under Condition 2--4, 1A, 1B and 2A, the bias of the estimator
$\widehat{\phi}_n(\lambda)$ given by (\ref{unifestn}) is
\begin{equation*}
\begin{split}
&\hskip-7pt
E[\widehat{\phi}_n(\lambda)-\phi(\lambda)]\\&=\left[-\frac{k_{q}}{2\pi}\int_{-\infty}^{\infty}|t|^{q}C(t)e^{-it\lambda}
dt\right](\rho_n b_n)^{q}+o\left((\rho_n b_n)^{q}\right)\\&~+
\left[-\frac{1}{2\pi}\int_{-\infty}^{\infty}|t|C(t)e^{-it\lambda}
dt\right]\left(\frac{\rho_n}{n}\right)+o\left(\frac{\rho_n}{n}\right)\\&~+
\left[\frac{A}{(2\pi)^{p}}\sum_{|l|>0}\frac{1}{|l|^p}\right]\frac{1}{(\rho_n)^{p}}+o\left(\frac{1}{(\rho_n)^{p}}\right),
\end{split}
\end{equation*}
i.e.,
$$E[\widehat{\phi}_n(\lambda)-\phi(\lambda)]=O\left((\rho_nb_n)^{q}\right)
+O\left(\frac{\rho_n}{n}\right)+O\left(\frac{1}{\rho_n^{p}}\right),
$$
uniformly in $\lambda$ over any closed and finite
interval.
\end{Thm}

\bigskip\noindent
{\bf Remark 7.} Condition 2A can be relaxed to the extent that the
characteristic exponent of the kernel $K(\cdot)$ is required to be
greater than or equal to the number $q$, for which Condition 1A is
assumed to hold. If it is strictly greater than $q$, then the term
$O\left((\rho_nb_n)^{q}\right)$ in the above theorem would have to
be replaced by $o\left((\rho_nb_n)^{q}\right)$. This follows from
equation (\ref{2arelax}) in the appendix and the fact that $k_q=0$
in this case. On the other hand, if a kernel with characteristic
exponent less than $q$ is used, then one does not fully utilize the
strength of the assumption on the smoothness of the spectral
density, implied by Condition 1A, and hence gets a slower rate of
convergence.

\subsection{Choice of $\rho_n$ and $b_n$ }\label{optrate}

From Theorem~3 and Theorem~2, it is observed that the bias and the
variance of the estimator $\widehat{\phi}_n(\lambda)$ converge to
zero at different rates. We set out to choose the sampling rate
$\rho_n$ and the window width $b_n$ in order to ensure that MSE of
$\widehat{\phi}_n(\lambda)$ converges to zero as fast as possible.
It would turn out that this happens when the squared bias and the
variance go to zero at the same rate.

\bigskip
\begin{Thm}
Under Conditions 2--5, 1A, 1B and 2A, the optimal rate of
convergence of the MSE of the estimator $\widehat{\phi}_n(\lambda)$
is given by
\begin{equation}
MSE[\widehat{\phi}_n(\lambda)]=O\left(n^{-\frac{2pq}{p+q+2pq}}\right), %
\label{msereg}
\end{equation}
which corresponds to the optimal choices
\begin{eqnarray}
\rho_n&=& P~ n^{\frac q{p+q+2pq}},\label{rhon} \\
\mbox{and }b_n&=&Q ~n^{-\frac{p+q}{p+q+2pq}},\label{bn}
\end{eqnarray}
for some positive constants $P$ and $Q$.
\end{Thm}

The above optimal rates of $\rho_n$ and $b_n$ lead to the following
corollaries to Theorems 2 and 3, respectively.

\bigskip\noindent
{\bf Corollary 1.} {\it Under the Conditions 1, 2, 5 and the choices
of $\rho_n$ and $b_n$ given by (\ref{rhon}--\ref{bn}), we have
\begin{equation}
\lim_{n\rightarrow\infty}n^{\frac{2pq}{p+q+2pq}}Var[\widehat{\phi}(\lambda)]
=\frac 1Q
(1+\delta_{0,\lambda})[\phi(\lambda)]^{2}\int_{-\infty}^{\infty}K^{2}(x)dx
,\qquad\label{phivar}%
\end{equation}
where $\delta_{0,\lambda}$ is equal to 1 if $\lambda=0$, and is 0
otherwise.}

\bigskip\noindent
{\bf Corollary 2.} {\it Under the Conditions 1A, 1B, 2, 2A and the
choices of $\rho_n$ and $b_n$ given by (\ref{rhon}--\ref{bn}), we
have
\begin{equation}
\begin{split}
&\hskip-15pt\lim_{n\rightarrow\infty}n^{\frac{pq}{p+q+2pq}}E[\widehat{\phi}_n(\lambda)-\phi(\lambda)]
\\=&-(PQ)^q \frac{k_{q}}{2\pi}\int_{-\infty}^{\infty}|t|^{q}C(t)e^{-it\lambda}
dt\!+\!\frac1{P^{p}}\frac{A}{(2\pi)^p}\sum_{|l|>0}\frac{1}{|l|^p},\label{phibias}
\end{split}
\end{equation}
where the constant $A$ is as in Condition 1B.}

\subsection{Comparison with Poisson Sampling }\label{poicomp}

Among the various schemes for sampling a continuous time stochastic
process at irregular intervals, the Poisson sampling proposed by
Silverman \cite{Silverman} is the simplest and most popular. Here, we compare
the asymptotic behaviour of the estimator $\widehat{\phi}_n(\cdot)$
with a similar estimator based on Poisson sampling.

Let $\{t_j\}_{j=0}^{n}$ be the sampling points from a Poisson
process with average sampling rate $\rho$. Masry \cite{Masrypoisson}
proved that, under Conditions 1, 2, 3 and 5, the estimator
$\widehat{\psi}_n(\cdot)$ defined as
\begin{equation*}
\begin{split}
&\hskip-7pt\widehat{\psi}_n(\lambda)\\ =&\frac{1}{\pi\rho
n}\sum_{i=1}^{n-1}\sum_{j=1}^{n-i}
X(t_j)X(t_{j+i})K(b_n(t_{j+i}\!-\!t_{j}))\cos(\lambda(t_{j+i}\!-\!t_{j})),
\end{split}
\end{equation*}
is consistent for $\phi(\lambda)$ for any choice of $\rho$.

Under the above conditions, the asymptotic variance of
$\widehat{\psi}_n(\lambda)$ satisfies
\begin{equation}\label{psivar}
\begin{split}
&\hskip-25pt\lim_{n\rightarrow\infty}(n b_n)Var[\widehat{\psi}_n(\lambda)]\\
=&\rho\left[\phi(\lambda)+\frac{C(0)}{2\pi\rho}\right]^{2}(1+\delta_{0,\lambda})\int_{-\infty}^{\infty}K^{2}(t)dt.
\end{split}
\end{equation}

For specifying the rate of convergence of the bias, Masry
\cite{Masrypoisson} assumed the following additional conditions.

\bigskip\noindent
{\bf Condition 1C.} $|t|^qC(t)$ is integrable for some positive
integer~$q$.

\bigskip\noindent
{\bf Condition 2B.} $K(\cdot)$ is $q$ times differentiable with
bounded derivatives, where $q$ is an integer for which Condition 1C
holds.

\bigskip
Note that Condition 1C is implied by Condition 1A with the same or
higher value of $q$ as is used here. Masry \cite{Masrypoisson}
showed that, under Conditions 1, 2, 3, 1C and 2B, the bias of the
estimator $\widehat{\psi}_n(\lambda)$ is given as
\begin{equation}\label{psibias}
\begin{split}
Bias[\widehat{\psi}_n(\lambda)]=&E[\widehat{\psi}_n(\lambda)]-\phi(\lambda)
\\=&\sum_{l=1}^{q-1}\frac{(i)^{l}K^{(l)}(0)b_n^{l}}{l!}\phi^{(l)}(\lambda)+O(b_n^{q})+O\left(\frac{1}{n}\right).
\end{split}
\end{equation}

It follows from (\ref{psibias}) that the rate of convergence of the
bias of $\widehat{\psi}_n(\lambda)$ is
$O\left(\max\left\{b_n^m,n^{-1}\right\}\right)$, where
$$m=\begin{cases} q&\mbox{if $K^{(l)}(0)=0$ for $1\le l<q$},\\
l_0&\mbox{if $K^{(l_0)}(0)\ne0$ and $K^{(l)}(0)=0$ for $1\!\le\!
l\!<\!l_0\!<\!q$}.\\\end{cases}$$

The fastest possible rate of convergence is
$O\left(\max\left\{b_n^q,n^{-1}\right\}\right)$, and this is
achieved when one uses a kernel, which further satisfies Condition
2A with the same or higher value of $q$ as is used here. In such a
case, we have
\begin{equation}
Bias[\widehat{\psi}_n(\lambda)]=O(b_n^{q})+O\left(\frac{1}{n}\right).
\label{psibias1}
\end{equation}
If the condition 2A holds with a higher value of $q$, then the term $O(b_n^{q})$ has to be replaced by $o(b_n^{q})$

Let us now assume that the kernel is chosen appropriately to ensure
(\ref{psibias1}). Note that the bias and the variance of
$\widehat{\psi}_n(\lambda)$ converge to zero at different rates. One
can choose the rate of convergence of the window width $b_n$ such
that the MSE converges as fast as possible. It turns out that if
$b_n=O(n^{-\alpha})$, then the optimal choice of $\alpha$ is
$\frac{1}{2q+1}$, in which case the squared bias and variance of
$\widehat{\psi}_n(\lambda)$ are both $O(n^{-\frac{2q}{2q+1}})$.

In summary, under Conditions 1, 1C, 2, 2A, 2B, 5 and
\begin{equation}
b_n=R~n^{-\frac1{2q+1}}, \label{bnP}
\end{equation}
the MSE of $\widehat{\psi}_n(\lambda)$ is
$$MSE[\widehat{\psi}_n(\lambda)]=O\left(n^{-\frac{2q}{2q+1}}\right).$$%
The rate of convergence of MSE of $\widehat{\phi}_n(\lambda)$ is
given in Theorem~4 under Conditions 1A, 1B, 2, 2A, 5 and
(\ref{rhon}--\ref{bn}). Both the results hold when $\rho_n$ and
$b_n$ for $\widehat{\phi}_n(\cdot)$ are chosen as in
(\ref{rhon}--\ref{bn}), $b_n$ for $\widehat{\psi}_n(\cdot)$ is
chosen as in (\ref{bnP}) and the following conditions hold
simultaneously: Condition 1A for some $q$ greater than 1 (which
implies Condition 1C for $[q]$ and Condition 1), Condition 2,
Condition 2A (for the same $q$ as in Condition 1A), condition 2B for
$[q]$, and Condition 5. Under this common set of conditions, we have
\begin{align*}
MSE[\widehat{\psi}_n(\lambda)]&=O\left(n^{-\frac{2[q]}{2[q]+1}}\right),\\
MSE[\widehat{\phi}_n(\lambda)]&=O\left(n^{-\frac{2q}{2q+1+q/p}}\right).
\end{align*}

When $q$ is an integer, i.e., $[q]=q$, the rate of convergence of
the MSE of $\widehat{\phi}_n(\lambda)$ is slower than that of
$\widehat{\psi}_n(\lambda)$. The two rates are comparable if $p$ is
much larger than $q$. When $q$ is not an integer, the MSE of
$\widehat{\phi}_n(\lambda)$ converges faster when $p>q[q]/(q-[q])$,
and in particular when $p$ is very large. As we have indicated in
Remark~6, for every fixed $q$, one can make $p$ suitably large
through low pass filtering.

If the rates of convergence are comparable, the constants associated
with these rates become important. We will compare the constants of
the asymptotic bias of $\widehat{\phi}_n(\lambda)$ and
$\widehat{\psi}_n(\lambda)$ as well as the constants of their
asymptotic variance separately, assuming that $q$ is an integer.

Under Conditions 1, 2, 5 and (\ref{bnP}), we have
\begin{equation}\label{psivar1}
\begin{split}
&\hskip-22pt\lim_{n\rightarrow\infty}n^{\frac{2q}{2q+1}}Var[\widehat{\psi}_n(\lambda)]\\
=&\frac1R\rho\left[\phi(\lambda)+\frac{C(0)}{2\pi\rho}\right]^{2}(1+\delta_{0,\lambda})\int_{-\infty}^{\infty}K^{2}(x)dx.
\end{split}
\end{equation}
On the other hand, we have from Corollary 1 that under the
Conditions 1, 2, 5 and (\ref{rhon}--\ref{bn}),
\begin{equation*}
\begin{split}
&\hskip-35pt\lim_{n\rightarrow\infty}n^{\frac{2q}{2q+1+q/p}}Var[\widehat{\phi}_n(\lambda)]\\
=&\frac 1Q
(1+\delta_{0,\lambda})[\phi(\lambda)]^{2}\!\int_{-\infty}^{\infty}K^{2}(x)dx.%\hskip25pt\mbox{}
\end{split}
\end{equation*}

The ratio of the constants for the asymptotic variances of
$\widehat{\psi}_n(\lambda)$ and $\widehat{\phi}_n(\lambda)$ is
\begin{equation}
\frac QR\rho\left[1+\frac{C(0)}{2\pi\rho\phi(\lambda)}\right]^{2}.
\label{psivaratio}
\end{equation}
This ratio depends on the Poisson sampling rate $\rho$ and the true
value of the power spectral density $\phi(\lambda)$. This ratio can
be much larger than 1, particularly for larger values of $\lambda$.
In fact, even if $\rho$ is chosen to minimize this ratio for a given
value of $\phi(\lambda)$ (though this is not practically possible),
the minimum value happens to be $2Q C(0)/[R \pi\phi(\lambda)]$, which
can be arbitrarily large for large values of $\lambda$. Thus, the
variance of $\widehat{\psi}_n(\lambda)$ can generally be expected to
be larger than that of $\widehat{\phi}_n(\lambda)$.

We now turn to the comparison of the expressions for bias. Under
Conditions 1C, 2, 2A and 2B along with (\ref{bnP}), we have
\begin{equation}\label{psibias2}
\begin{split}
&\hskip-40pt\lim_{n\rightarrow\infty}n^{\frac{q}{2q+1}}[E[\widehat{\psi}_n(\lambda)]-\phi(\lambda)]\\
=&-R^qk_{q}\frac{1}{2\pi}\int_{-\infty}^{\infty}|t|^{q}C(t)e^{-it\lambda}
dt.
\end{split}
\end{equation}
On the other hand, we have from Corollary 2 that under the
Conditions 1A, 1B, 2, 2A and (\ref{rhon}--\ref{bn}),
\begin{eqnarray*}
\begin{split}
&\hskip-10pt\lim_{n\rightarrow\infty}n^{\frac{pq}{p+q+2pq}}E[\widehat{\phi}_n(\lambda)-\phi(\lambda)]
\\=&-(PQ)^q k_{q}\frac{1}{2\pi}\int_{-\infty}^{\infty}|t|^{q}C(t)e^{-it\lambda}
dt+\frac1{P^p}\frac{A}{(2\pi)^p}\sum_{|l|>0}\frac{1}{|l|^p}.
\end{split}
\end{eqnarray*}
The first term of the expression on the right hand side is
proportional to the expression on the right hand side of
(\ref{psibias2}). These terms are small for large values of $\lambda$, while the
second term of the expression on the right hand side of the above
inequality does not depend on $\lambda$. Even if the value of the
second term is small, it would make a difference for large values of
$\lambda$. Consequently, $\widehat{\psi}_n(\lambda)$ can generally
be expected to have a smaller bias than $\widehat{\phi}_n(\lambda)$.

In summary, even though both $\widehat{\phi}_n(\lambda)$ and
$\widehat{\psi}_n(\lambda)$ are consistent estimators under the
stated conditions, there is a trade-off between
$\widehat{\phi}_n(\lambda)$ and $\widehat{\psi}_n(\lambda)$ in terms
of bias and variance. There is no clear order between the constants
of the mean square errors of the two estimators.

In order to examine the validity of the asymptotic results and the
above comparisons for small samples, we turn to Monte Carlo
simulations, reported in the next section.

\section{Simulation study}\label{s4}
In this section, we shall present the simulation study of
performance of the spectral density estimators based on regular and
Poisson sampled data. We consider a continuous time
autoregressive (AR(4)) process having the spectral density
$$\phi(\lambda)=\frac{\sigma^{2}}
{2\pi}\cdot\frac1{(\lambda^2+\alpha_1^2)(\lambda^2+\alpha_2^2)(\lambda^2+\alpha_3^2)(\lambda^2+\alpha_4^2)},$$
where $\alpha_1=0.65$, $\alpha_2=0.75$, $\alpha_3=0.85$,
$\alpha_4=0.95$ and $\sigma=1$.

A process $\{X(t),~ -\infty<t<\infty\} $ having the above spectral
density can be written as \cite{Hoel}
$$X(t)=\int_{-\infty}^{t}h(t-s)dW(s),$$
where the impulse function $h(\cdot)$ is given by
$$h(t)=\sum_{i=1}^{4}c_i e^{-\alpha_i t}1_{[0,\infty)}(t)
$$
and the constants $c_i$, $i=1,\ldots,4$, are the solution of the
following system of linear equations:
$$h(0)=h^{(1)}(0)=h^{(2)}(0)=0;\quad h^{(3)}(0)=1,$$
$h^{j}(0)$ being the $j$th derivative of $h(\cdot)$ evaluated at
$0$.

In view of the above representation, we simulate the process $\{X_0(t)\}$ given by
$$X_0(t)=\int_{0}^{t}h(t-s)dW(s).$$
The process $\{X_0(t),~ 0<t<\infty\}$ is not a stationary process.
However, as $t$ becomes large, the variance of the difference
between the processes $\{X_0(t)\}$ and $\{X(t)\}$ becomes small. We
find out the value of $t$, say $t_0$, such that
$Var(X_0(t)-X(t))<10^{-9}$, and consider the path of the simulated
process $\{X_0(t)\}$ from $t_0$ onwards.

For estimation, we assume that the underlying power spectral density
satisfies Condition 1A with $q=2$. According we use the Hanning
Kernel
$$K(x)=\frac{1}{2}(1+\cos(\pi x))1_{[-1,1]}(x),$$
which has characteristic exponent~2.

\subsection{Finite sample performance of the estimator $\widehat{\phi}_n(\lambda)$}
Here, we consider the performance of $\widehat{\phi}_n(\lambda)$
over the frequency range $[0,\pi/2]$. We used the optimal choice of
sampling rate developed in Section~\ref{optrate} to generate
regularly spaced samples of the process for sample sizes
$n=100,~1000~\mbox{and}~10000$. We assume Condition 1A with $q=2$
and Condition 1B with $p=8$ (both of which actually hold for the
underlying power spectral density). For the above choices, the
optimal powers of $n$ for the sampling rate and the window width are
$\rho_n\propto n^{1/21}$ and $b_n\propto n^{-5/21}$. We choose
$\rho_n= n^{1/21}$ and $b_n=\frac{1}{4}n^{-5/21}$.

Figure 1 shows the average of the estimated power spectral density
computed from 500 simulation runs, the empirically observed bias and
variance, together with the true power spectral density and the
theoretical (asymptotic) bias and variance, respectively, for the
three samples sizes.
\begin{figure*}
\noindent\mbox{}\hskip-15pt
\includegraphics[width=8in,height=3in]{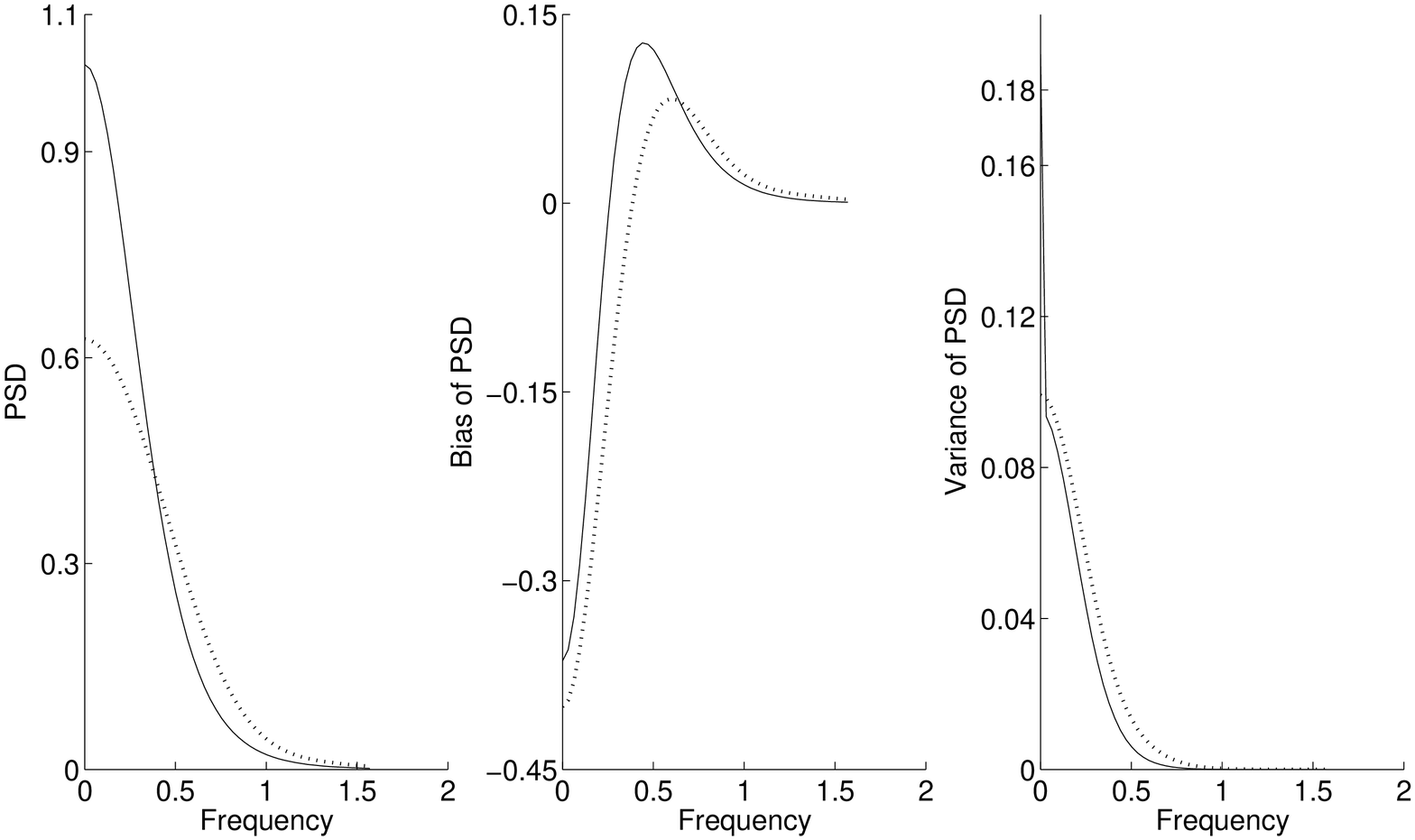}\\\mbox{}\hskip-15pt
\includegraphics[width=8in,height=3in]{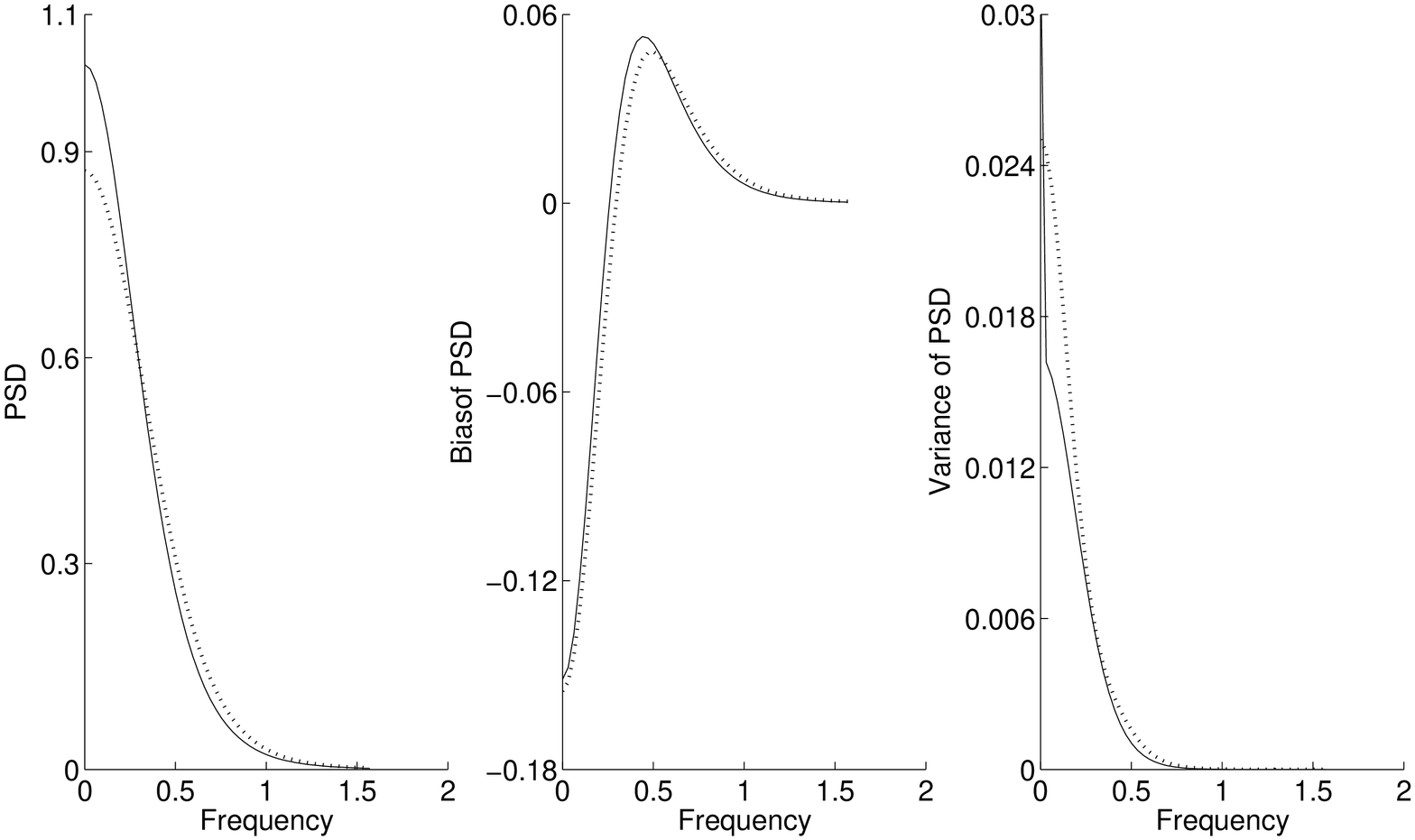}\\\mbox{}\hskip-15pt
\includegraphics[width=8in,height=3in]{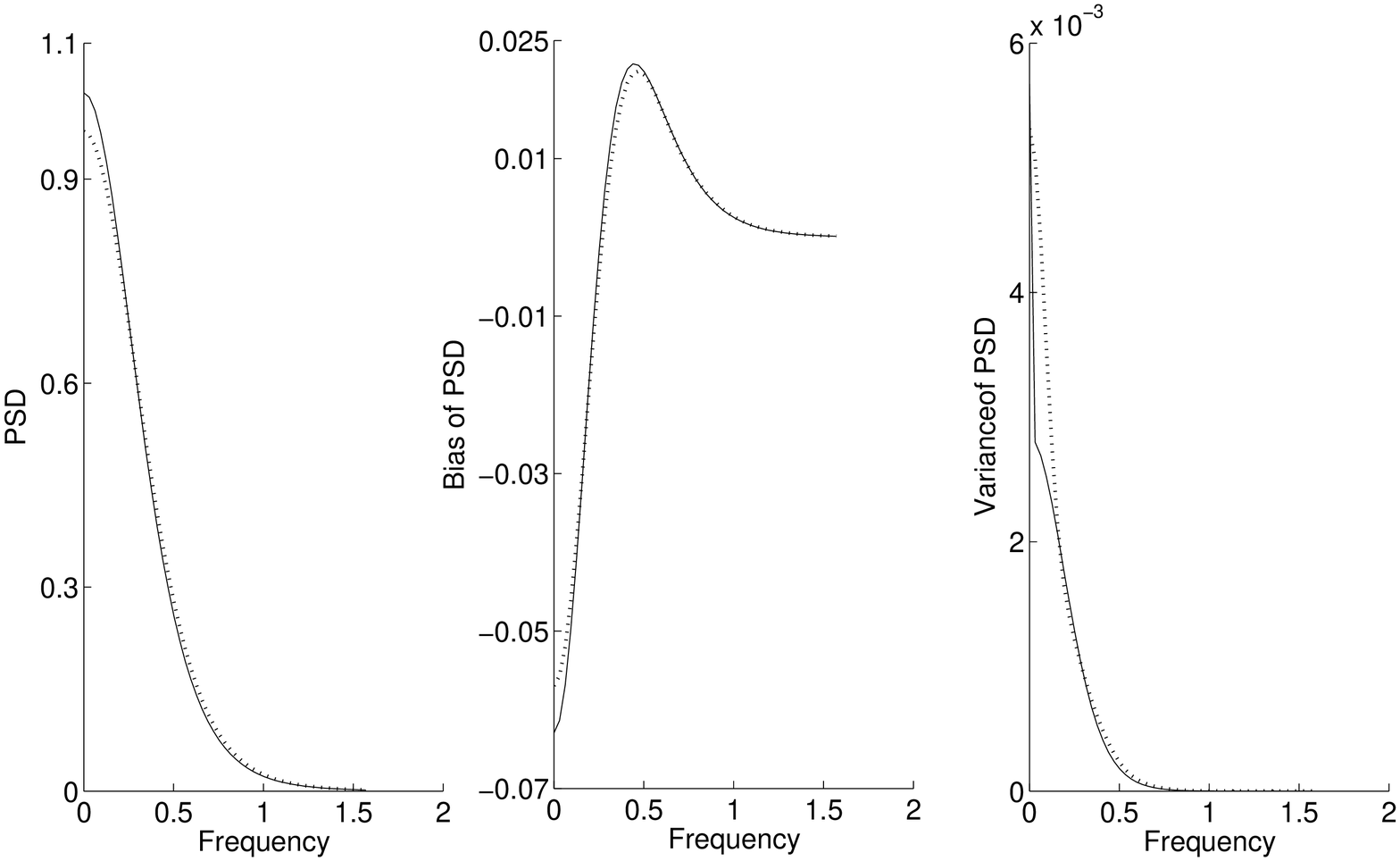}
\caption{The average estimated power spectral density
$\widehat{\phi}_n(\cdot)$ (left column), the bias (middle column)
and the variance (right column) for sample sizes 100 (top row), 1000
(middle row) and 10000 (bottom row). The solid and the dotted lines
correspond to theoretical (asymptotic) and empirical values,
respectively.}
\end{figure*}

From these figures, it can be observed that as the sample size goes
from 100 to 10000, the empirical values of bias and variance get
closer to the asymptotic results. Moreover, the theoretical
(asymptotic) computations are quite comparable to the empirical
values, even for sample size 100.

\subsection{Finite sample comparison of $\widehat{\phi}_n(\lambda)$
with Poisson sampled estimator $\widehat{\psi}_n(\lambda)$
}\label{4b}

We generate Poisson sampled data with the average sampling rate
$\rho=1$ for sample sizes $n=100$, 1000 and 10000, and compute the
estimator $\widehat{\psi}_n(\lambda)$ on $[0,\pi/2]$. Here, the
optimal power of $n$ for the window width is $b_n \propto n^{-1/5}$.
We use $b_n=\frac{1}{4}n^{-1/5}$.

Figure 2 shows the empirical bias, variance and MSE of the
estimators $\widehat{\phi}_n(\lambda)$ and
$\widehat{\psi}_n(\lambda)$ computed from 500 simulation runs, as a
function of the frequency, for sample sizes 100, 1000 and 10000.

\begin{figure*}
\noindent\mbox{}\hskip-15pt
\includegraphics[width=8in,height=3in]{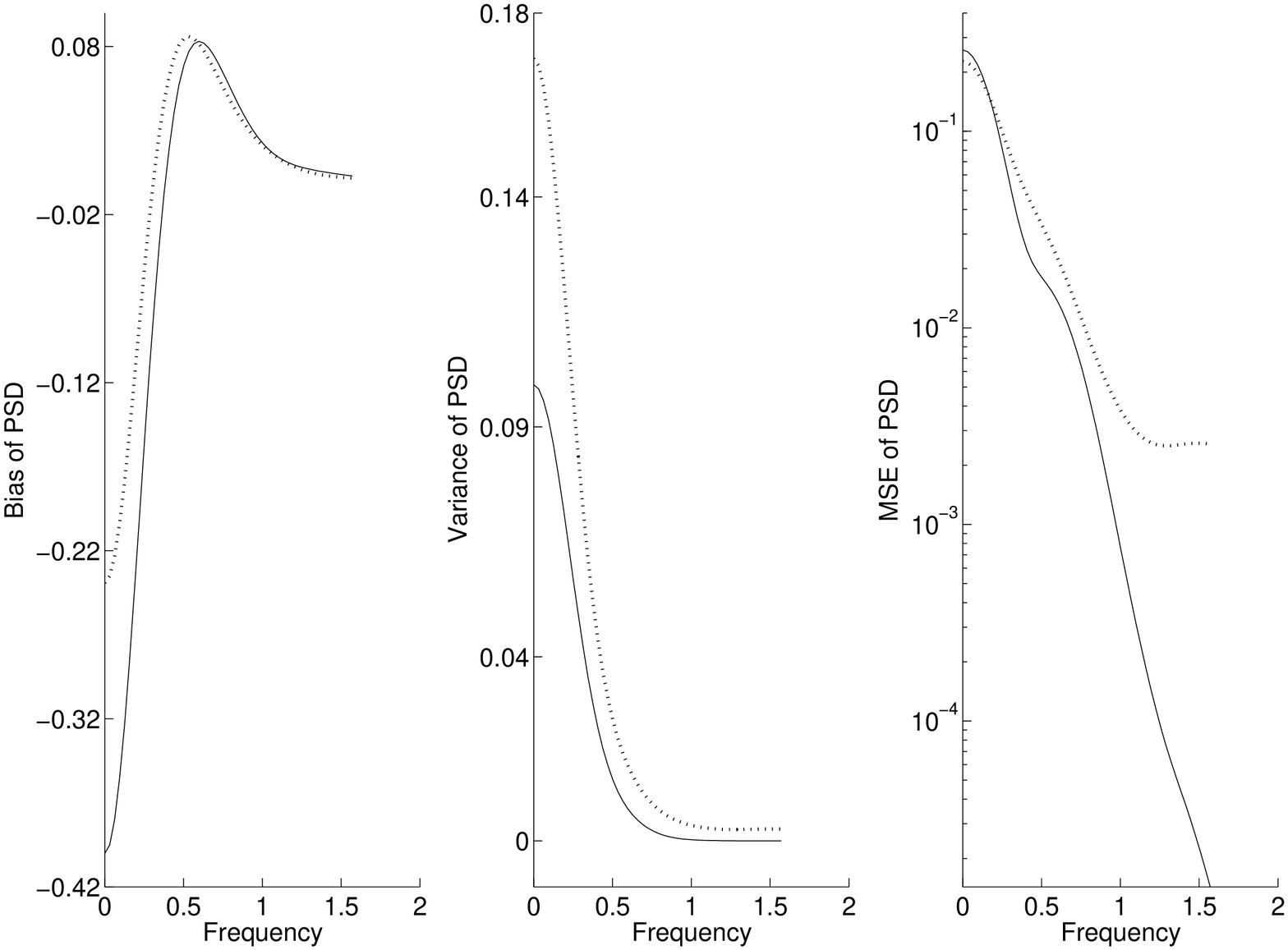}\\\mbox{}\hskip-15pt
\includegraphics[width=8in,height=3in]{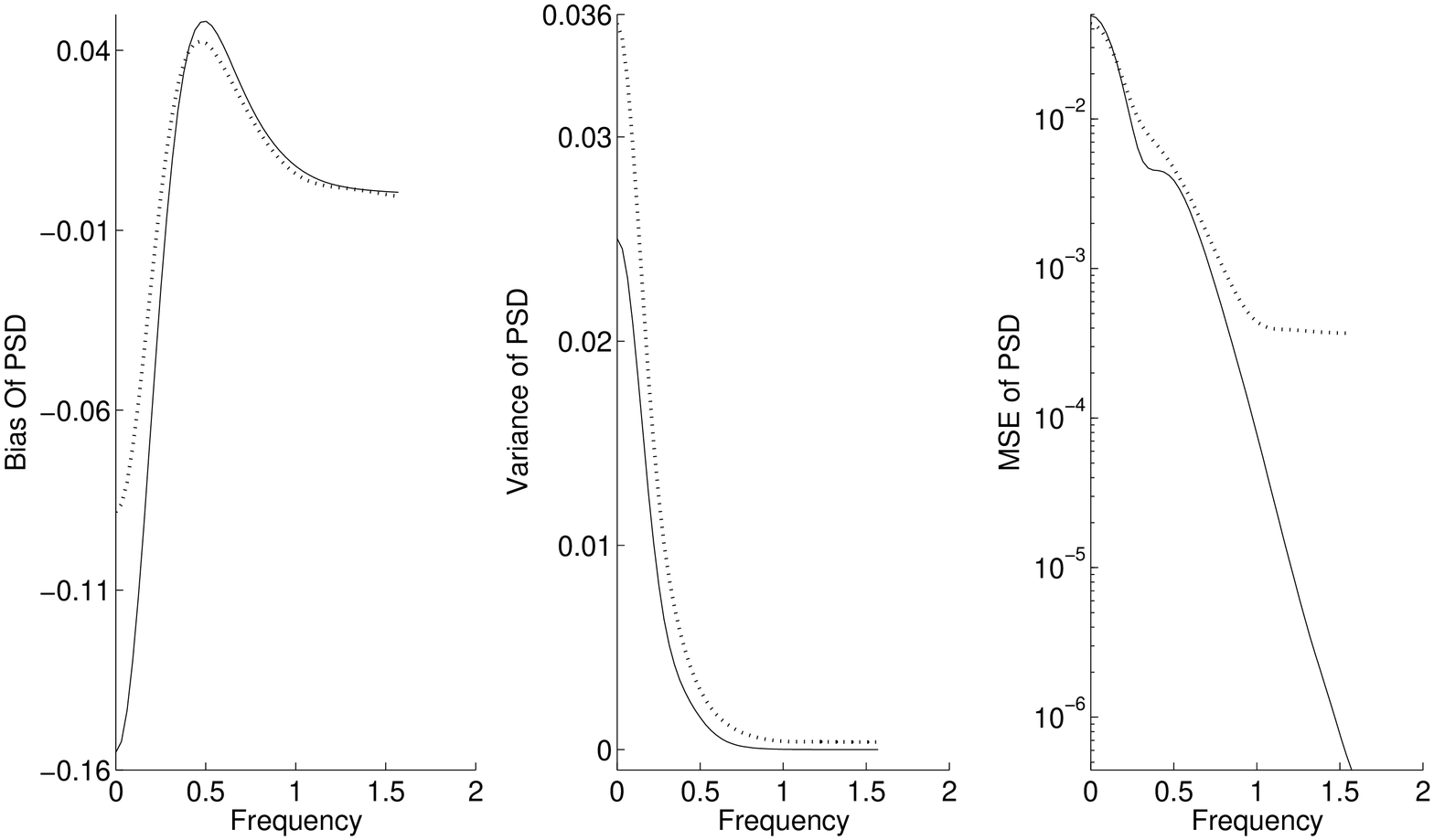}\\\mbox{}\hskip-15pt
\includegraphics[width=8in,height=3in]{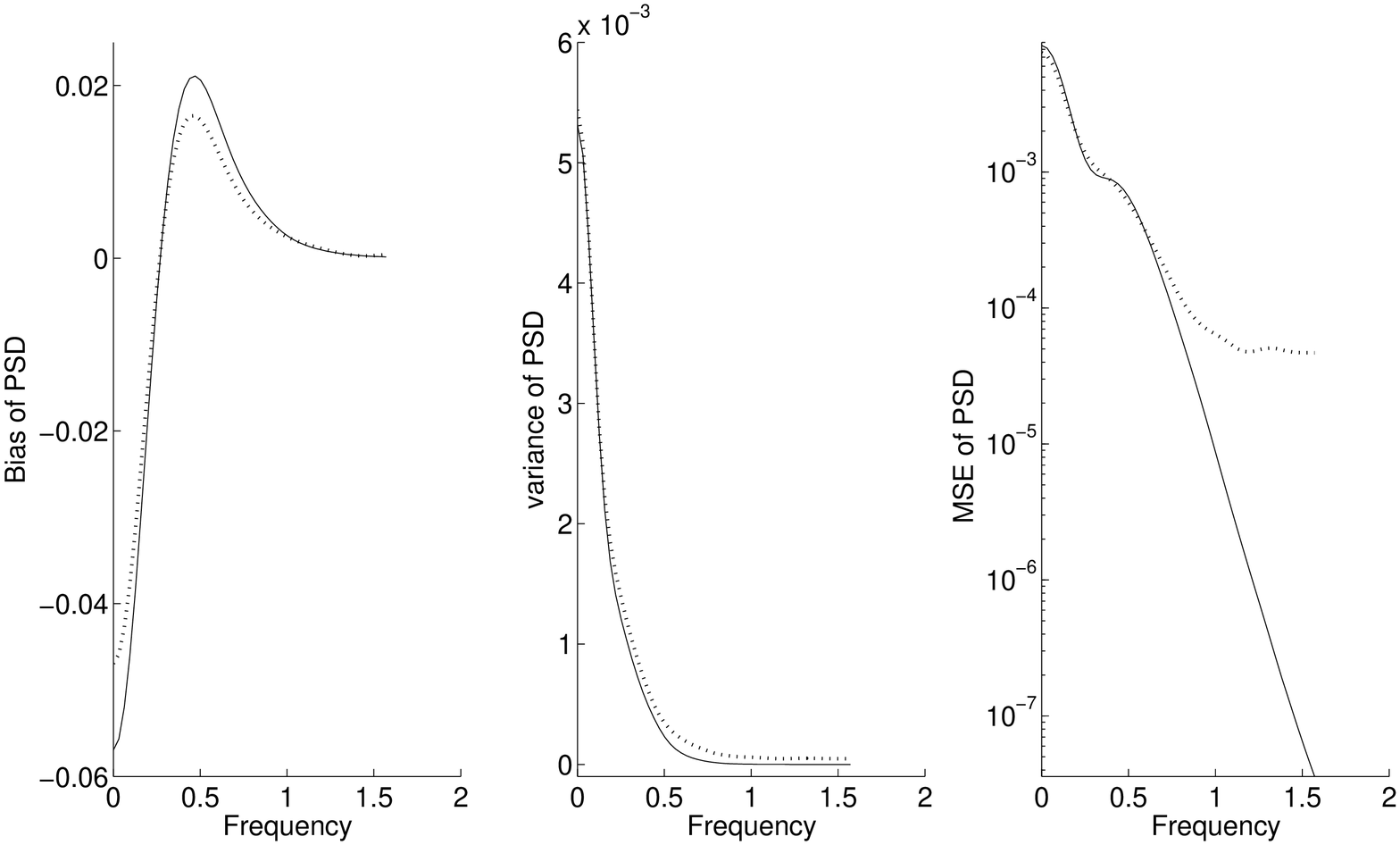}
\caption{The estimated bias of the power spectral density
$\widehat{\phi}_n(\cdot)$ and $\widehat{\psi}_n(\cdot)$ (left
column), the variance (middle column) and the MSE (in log scale,
right column) for sample sizes 100 (top row), 1000 (middle row) and
10000 (bottom row). The solid and the dotted lines correspond to
regular and Poisson samplings, respectively.}
\end{figure*}

From these figures, it can be observed that the bias of
$\widehat{\psi}_n(\cdot)$ is generally less than that of
$\widehat{\phi}_n(\cdot)$ while the variance of
$\widehat{\psi}_n(\cdot)$ is larger than that of
$\widehat{\phi}_n(\cdot)$. The differences diminish with larger
sample size. These patterns are in accordance with the large sample
comparisons made in Section~\ref{poicomp}. The MSE of
$\widehat{\psi}_n(\cdot)$ is larger than that of the
$\widehat{\phi}_n(\cdot)$ for larger frequencies. The MSE is plotted
in log-scale in order to highlight the fact that this quantity, in
the case of $\widehat{\psi}_n(\cdot)$, levels off to a constant
value for larger frequencies, while in the case of
$\widehat{\phi}_n(\cdot)$, it continues to decline. This difference
in behaviour is in accordance with the variance expressions given in
(\ref{phivar}) and (\ref{psivar1}).

\section{Discussion}\label{s5}

In this paper, we have shown that the smoothed periodogram based on
regularly spaced samples of a continuous time stationary stochastic
process is consistent, under certain conditions, provided that the
sampling rate increases appropriately as the sample size goes to
infinity. We have also shown that, under the conditions used in the
proofs, the estimators based on uniformly and non-uniformly spaced
samples have about the same rates of convergence. Thus, our results
remove a widely perceived theoretical deficiency of a popular
spectral estimator based on regular sampling.

It has been a common experience, both theoretically and
empirically\cite{Robert}, that the smoothed periodogram estimator
(\ref{unifest}) of a non-bandlimited power spectral density has less
variance and more bias compared to the corresponding estimator
$\widehat{\psi}_n(\cdot)$ based on Poisson sampling. What the
results of Section~\ref{s3} show is that, even though the new
asymptotic results presented in this paper establish consistency of
the smoothed periodogram $\widehat{\phi}_n(\cdot)$ and the rates of
convergence of the estimators $\widehat{\phi}_n(\cdot)$ and
$\widehat{\psi}_n(\cdot)$ are comparable, the constants for the
first order approximations of the bias and variance of the two
estimators exhibit the same type of trade off, i.e., the constant
for the bias term is larger in the case of
$\widehat{\phi}_n(\cdot)$, and the constant for the variance term is
larger in the case of $\widehat{\psi}_n(\cdot)$.

The new asymptotic calculations provide a theoretical justification
of using the smoothed periodogram with common sense, even if the
underlying power spectral density is not bandlimited. This common
sense approach consists of appropriate filtering of the continuous
time process followed by sampling at a suitably uniform rate.
Remark~6 and Theorem~4 give guidelines for choosing a suitable
filter and an appropriate sampling rate, respectively, which may be
useful for practitioners.

The simulation results reported in Section~\ref{s4} illustrate how
one can choose an appropriate sampling rate for estimating the power
spectral density, and obtain results in line with the theoretical
results. Even though the underlying spectral density in this example
is not band-limited, the estimator $\widehat\phi_n(\cdot)$ (based on
uniformly spaced samples) is found to have smaller MSE than
$\widehat\psi_n(\cdot)$ (based on Poisson samples) for larger
frequencies. The reverse order holds for smaller frequencies. This
shows that there is no clear dominance of one kind of sampling over
another. This finding for finite samples complements our asymptotic
results.

Our results do not take anything away from the vast literature on
spectrum estimation through irregularly sampled data. These methods
may be quite appropriate when one does not have control over the
sampling mechanism, when irregular sampling is logistically feasible
and methodologically not limited, or when regular sampling have to
be avoided for a specific reason (other than its perceived
inconsistency). Further, an irregular sampling scheme such as
Poisson sampling can be used where an estimator based on it is
expected, either through theoretical analysis or through simulation
studies, to have smaller MSE than the corresponding estimator based
on regular sampling scheme.

We have shown in Section~\ref{poicomp} how our theoretical results
can be used to compare uniform and Poisson sampling schemes. The
results compiled there may be used to make further comparisons under
different constraints. For example, if there is a limit to the
maximum average sampling rate and/or the maximum sample size, one
may make an optimal choice of the window width for fixed values of
these two parameters, and then determine the corresponding MSE. In
the case of $\widehat{\phi}_n(\cdot)$, the optimal choice of the
window width (for given sample size and sampling rate) is given by
(\ref{interm}), while the choice in the case of
$\widehat{\psi}_n(\cdot)$ can be derived similarly from
(\ref{psivar}) and (\ref{psibias})\cite{Masrypoisson}. The best
rates and constants achievable under the two sampling schemes, under
appropriate constraints, may then be used to make a suitable choice
of the sampling scheme. However, if there is a hard restriction on
the minimum separation between two successive samples (rather than a
restriction on the {\it average} sampling rate), then one cannot use
Poisson sampling at all. In such cases, irregular sampling may be
done according to a renewal process, with the inter-sample distance
having a restricted probability distribution. Such a sampling scheme
would not satisfy the sufficient condition for alias-free sampling
given in Theorem 1 of \cite{Masrynew}. Thus, one may have to look
further in search of a suitable estimator based on non-uniform
sampling under such a restriction.

The proven consistency of the smoothed periodogram opens up the
possibility of establishing consistency of {\it parametric}
estimators of the power spectral density of a continuous time
process based on regularly spaced samples, by allowing the sampling
rate together with the sample size to go to infinity. Such
asymptotic calculations may potentially be used to justify and/or
fine-tune multi-resolution methods of spectrum
estimation\cite{Eldar}.

\appendix\setcounter{section}{1}
We denote by $K_1(\cdot)$ a function that bounds the covariance
averaging kernel $K(\cdot)$ as in Condition~2. Further, we denote
$K_1(0)$ by $M$.

\bigskip\noindent
{\bf Proof of Theorem 1.}\ \ We shall show that the bias of the estimator $\widehat{\phi}_n(\lambda)$ given by (\ref{unifestn}) converges to $0$ uniformly over $[\lambda_l,\lambda_u]$ for any $\lambda_l$, $\lambda_u$ such that $\lambda_l<\lambda_u$. In order to compute the bias, we
evaluate $E[\widehat{\gamma}_n(v)]$:
\begin{align}
E[\widehat{\gamma}_n(v)]&=E\left[\frac{1}{n}\sum_{j=1}^{n-|v|}X
\left(\frac{j}{\rho_n}\right)X\left(\frac{j+|v|}{\rho_n}\right)\right]\notag\\
&=\left(1-\frac{|v|}{n}\right)C\left(\frac{v}{\rho_n}\right).
\label{expg}
\end{align}
Therefore, we have
\begin{align*}
&\hskip-7pt
E[\widehat{\phi}_n(\lambda)]\\
=&\frac{1}{2\pi\rho_n}\sum_{|v|<n}\!\left(1\!-\!\frac{|v|}{n}\right)C\left(\frac{v}{\rho_n}\right)
K(b_n
v)e^{\frac{-iv\lambda}{\rho_n}}1_{[-\pi\rho_n,\pi\rho_n]}(\lambda).
\end{align*}
Consider the simple function $S_n(\cdot)$, defined over $[\lambda_l,\lambda_u]\times(-\infty,\infty)$, by
\begin{align*}
S_n(\lambda,t)=&\frac{1}{2\pi}\sum_{|v|<n}\left(1-\frac{|v|}{n}\right)C\left(\frac{v}{\rho_n}\right)K\left(b_n
v\right)e^{\frac{-iv\lambda}{\rho_n}}\\&\times1_{[-\pi\rho_n,\pi\rho_n]}(\lambda)1_{\left(\frac{v-1}{\rho_n},\frac{v}{\rho_n}\right]}(t).
\end{align*}
Observe that $\int_{-\infty}^{\infty} S_n(\lambda,t)dt=E[\widehat{\phi}_n(\lambda)]$.\\
Define the function $S(\cdot)$, over $[\lambda_l,\lambda_u]\times(-\infty,\infty)$, by
$$ \hskip-80pt S(\lambda,t)=\frac{1}{2\pi}C(t)e^{-it\lambda}.$$
Observe that $\int_{-\infty}^{\infty} S(\lambda,t)dt=\phi(\lambda)$ which is continuous.

For any $t \in (-\infty,\infty)$, let $v_n(t)$ be the smallest
integer greater than or equal to $\rho_nt$. Note that the interval
$(\frac{v_{n-1}(t)}{\rho_n},\frac{v_n(t)}{\rho_n}]$ contains the
point $t$  and $\lim_{n\rightarrow\infty}\frac{v_n(t)}{\rho_n}= t$.
For sufficiently large $n$, we have from Conditions 3 and 4,
\begin{align*}
S_n(\lambda,t)=&\frac{1}{2\pi}\left(1-\frac{|v_n(t)|}{\rho_n}\cdot\frac{\rho_n}{n}\right)
C\left(\frac{v_n(t)}{\rho_n}\right)K\left(b_n \rho_n\frac{v_n(t)}{\rho_n}\right)\\
&\times
e^{\frac{-iv_n(t)\lambda}{\rho_n}}1_{[-\pi\rho_n,\pi\rho_n]}(\lambda).
\end{align*}
Proving the uniform convergence of $Bias[\widehat{\phi}_n(\lambda)]$
over finite interval $[\lambda_l,\,\lambda_u]$ amounts
to proving
$$\lim_{n\rightarrow\infty}\int_{-\infty}^{\infty} S_n(\lambda,t)dt=\int_{-\infty}^{\infty} S(\lambda,t)dt,$$
uniformly over  $[\lambda_l,\,\lambda_u]$. By virtue of the
continuity of the limiting function, this in turn is equivalent to
proving that $\int_{-\infty}^{\infty} S_n(\lambda,t)dt$ converges
continuously over this interval \cite{Resnick}, i.e., for any
sequence $\lambda_n \rightarrow \lambda$,
$$\lim_{n\rightarrow\infty}\int_{-\infty}^{\infty} S_n(\lambda_n,t)dt = \int_{-\infty}^{\infty} S(\lambda,t)dt$$
where $\lambda_n,\lambda \in  [\lambda_l,\,\lambda_u]$.

By continuity of the function $S_n(\lambda,t)$ with respect to $t$
and $\lambda$, we have from Conditions 3 and 4, for any fixed $t$,
$$\lim_{n\rightarrow \infty}|S_n(\lambda_n,t)-S(\lambda,t)|=0.$$
Note that from Conditions 1 and 2, we have the dominance
\begin{align*}
|S_n(\lambda_n,t)|\le&
M\sum_{|v|<n}\left|C\left(\frac{v}{\rho_n}\right)\right|1_{\left(\frac{v-1}{\rho_n},\frac{v}{\rho_n}\right]}(t)\le
M h_0(t).
\end{align*}
where $h_0(\cdot)$ is the function described in Condition~1.
Thus, by applying the dominated convergence theorem (DCT), we have
$$\lim_{n\rightarrow \infty}\int_{-\infty}^{\infty}
S_n(\lambda_n,t)dt= \int_{-\infty}^{\infty} S(\lambda,t)dt.$$

Hence, $E[\widehat{\phi}_n(\lambda)]\rightarrow \phi(\lambda)$
uniformly on $[\lambda_l,\,\lambda_u]$. \hfill$\Box$

\bigskip
\noindent {\bf Proof of Theorem 2.}\ \ The estimator
$\widehat{\phi}_n(\lambda)$, given by (\ref{unifestn}), can be
written as
$$\widehat{\phi}_{n}(\lambda)=\frac{1}{2\pi\rho_n}\widehat{\gamma}_n(0)+\frac{1}{\pi \rho_n}\sum_{v=1}^{n-1}\widehat{\gamma}_n(v)K(b_n v)\cos\left(\frac{\lambda v}{\rho_n}\right).$$
Therefore,
\begin{equation}
Var[\widehat{\phi}_{n}(\lambda)]=I_1+2I_2(\lambda)+I_3(\lambda),
\end{equation}
where
\begin{IEEEeqnarray*}{rCl}
I_1&=&\frac{1}{(2\pi)^{2}\rho_n^{2}}Var[\widehat{\gamma}_n(0)],\\
I_2(\lambda)&=&Cov\left[\frac{1}{2\pi
\rho_n}\widehat{\gamma}_n(0),\frac{1}{\pi
\rho_n}\sum_{v=1}^{n-1}\widehat{\gamma}_n(v)K(b_n
v)\cos\!\left(\frac{\lambda v}{\rho_n}\right)\!\right]\!,\\
I_3(\lambda)&=&Var\left[\frac{1}{\pi
\rho_n}\sum_{v=1}^{n-1}\widehat{\gamma}_n(v)K(b_n
v)\cos\left(\frac{\lambda v}{\rho_n}\right)\right]\!.
\end{IEEEeqnarray*}
Before we consider the convergence of the above three terms, we
simplify the computation of
$Cov(\widehat{\gamma}_n(v_{1}),\widehat{\gamma}_n(v_{2}))$ for non
negative $v_1$ and $v_2$. Note from Condition~5 that
\begin{equation}\label{exp_cros}
\begin{split}
&\hskip-7pt E[\widehat{\gamma}_n(v_{1})\widehat{\gamma}_n(v_{2})]\\
=&E\left[\frac{1}{n^{2}}\sum_{j=1}^{n-v_{1}}\sum_{l=1}^{n-v_{2}}X\!\!\left(\frac{j}{\rho_n}\right)
X\!\!\left(\frac{j\!+\!v_{1}}{\rho_n}\right)X\!\!\left(\frac{l}{\rho_n}\right)
X\!\!\left(\frac{l\!+\!v_2}{\rho_n}\right)\right]\\
=&\frac{1}{n^{2}}\sum_{j=1}^{n-v_{1}}\sum_{l=1}^{n-v_{2}}
\left[C\!\!\left(\frac{v_{1}}{\rho_n}\!\right)\!C\!\!\left(\frac{v_{2}}{\rho_n}\!\right)
+C\!\!\left(\!\frac{j\!-\!l\!+\!v_{1}}{\rho_n}\!\right)\!C\!\!\left(\!\frac{j\!-\!l\!-\!v_{2}}{\rho_n}\!\right)\right.\\
&\left.~
+C\!\!\left(\!\frac{j\!-\!l}{\rho_n}\right)\!C\!\!\left(\!\frac{j\!-\!l\!+\!v_{1}\!-\!v_{2}}{\rho_n}\right)
+Q\!\!\left(\frac{v_1}{\rho_n},\frac{l\!-\!j}{\rho_n},\frac{l\!-\!j\!+\!v_2}{\rho_n}\right)\right]\\
=&\left(1-\frac{v_{1}}{n}\right)C\!\!\left(\frac{v_{1}}{\rho_n}\right)
\left(1-\frac{v_{2}}{n}\right)C\!\!\left(\frac{v_{2}}{\rho_n}\right)\\
&~+\frac{1}{n}\sum_{|u|<n}U_{n}(u,v_{1},v_{2})~C\!\!\left(\frac{u+v_{1}}{\rho_n}\right)C\!\!\left(\frac{u-v_{2}}{\rho_n}\right)\\
&~+\frac{1}{n}\sum_{|u|<n}U_{n}(u,v_{1},v_{2})~C\!\!\left(\frac{u}{\rho_n}\right)C\!\!\left(\frac{u+v_{1}-v_{2}}{\rho_n}\right)
\\
&~+\frac{1}{n}\sum_{|u|<n}U_{n}(u,v_{1},v_{2})~Q\!\!\left(\frac{v_1}{\rho_n},\frac{-u}{\rho_n},\frac{-u+v_2}{\rho_n}\right),
\end{split}
\end{equation}
where $U_{n}(u,v_{1},v_{2})$ is a function with values between 0 and
1 defined as follows
\begin{equation}\label{desofu}
\begin{split}
&\hskip-7pt U_{n}(u,v_{1},v_{2})\\
=&
\begin{cases}
0,  &u \leq -n+v_{2},\\
1-\frac{v_{2}\!-\!u}{n} &-n\!+\!v_{2}\!<\! u\! <\! \min(0,v_{2}\!-\!v_{1}),\\
1-\frac{\max(v_{1},v_{2})}{n} &\min(0,v_{2}\!-\!v_{1})\leq u \leq \max(0,v_{2}\!-\!v_{1}),\\
1-\frac{u\!+\!v_{1}}{n} &\max(0,v_{2}\!-\!v_{1}) < u < n\!-\!v_{1},\\
0 &u\geq n\!-\!v_{1}.
\end{cases}
\end{split}
\end{equation}
Therefore, by using (\ref{expg}) and (\ref{exp_cros}), we have
\begin{equation}\label{desofcov}
\begin{split}
&\hskip-10ptCov(\widehat{\gamma}_n(v_{1}),\widehat{\gamma}_n(v_{2}))\\
=& E[\widehat{\gamma}_n(v_{1})\widehat{\gamma}_n(v_{2})]-E[\widehat{\gamma}_n(v_{1})]E[\widehat{\gamma}_n(v_{2})]\\
=&\frac{1}{n}\sum_{|u|<n}U_{n}(u,v_{1},v_{2})
C\left(\frac{u+v_{1}}{\rho_n}\right)C\left(\frac{u-v_{2}}{\rho_n}\right)\\
&+\frac{1}{n}\sum_{|u|<n}U_{n}(u,v_{1},v_{2})C\left(\frac{u}{\rho_n}\right)C\left(\frac{u+v_{1}-v_{2}}{\rho_n}\right)
\\&+
\frac{1}{n}\sum_{|u|<n}U_{n}(u,v_{1},v_{2})Q\left(\frac{v_1}{\rho_n},\frac{-u}{\rho_n},\frac{-u+v_2}{\rho_n}\right).
\end{split}
\end{equation}
We now use this simplified form of
$Cov(\widehat{\gamma}_n(v_{1}),\widehat{\gamma}_n(v_{2}))$ to
establish the convergence of the three terms, $I_1$, $I_2(\lambda)$
and $I_3(\lambda)$.

Using (\ref{desofcov}) and Condition 5, $I_1$ is given as
\begin{align*}
I_1
=&\frac{1}{(2\pi)^{2}\rho_n^{2}}\left[\frac{2}{n}\sum_{|u|<n}U_{n}(u,0,0)C\left(\frac{u}{\rho_n}\right)^{2}\right.
\\&\left.+\frac{1}{n}\sum_{|u|<n}U_{n}(u,0,0)Q\left(0,\frac{-u}{\rho_n},\frac{-u}{\rho_n}\right)\right] \\
\le&\frac{2C(0)}{(2\pi)^{2}n\rho_n}\sum_{|u|<n}\!\left|C\!\!\left(\frac{u}{\rho_n}\!\right)\!\right|\frac{1}{\rho_n}
+\frac{g_1(0)g_2(0)}{(2\pi)^{2}n\rho_n}\!\sum_{|u|<n}\!\!g_3\!\!\left(\frac{u}{\rho_n}\!\right)\!\!\frac{1}{\rho_n}.
\end{align*}
As in Theorem 1, we can view
$\sum_{|u|<n}\!\left|C\!\!\left(\frac{u}{\rho_n}\!\right)\!\right|\frac{1}{\rho_n}$
as the integral of the function $s_n(\cdot)$ defined by
$$s_n(t)=\sum_{|i|<n}\left|C\left(\frac{i}{\rho_n}\right)\right|1_{(\frac{i-1}{\rho_n},\frac{i}{\rho_n}]}(t).$$
Since $s_n(t)\!\rightarrow\! |C(t)|$, and $s_n(t)\!\le\! h_0(t)$ holds
from Condition~1, we get
$\lim_{n\rightarrow\infty}\sum_{|u|<n}\left|C\!\!\left(\frac{u}{\rho_n}\right)\right|\frac{1}{\rho_n}\!=\!
\int_{-\infty}^{\infty}|C(t)|dt$ by applying the DCT as in
Theorem~1, under Condition~4. A similar argument, together with
Conditions 4 and 5, ensures that
$\lim_{n\rightarrow\infty}\sum_{|u|<n}g_3\left(\frac{u}{\rho_n}\right)\frac{1}{\rho_n}=
\int_{-\infty}^{\infty}|g_3(u)|du$. Both the limiting integrals are
finite. So $nb_nI_1\rightarrow 0$ as $n\rightarrow\infty$.

\medskip
Using (\ref{desofcov}) and Condition 5, the term $I_2(\lambda)$ is
given as
\begin{equation*}
\begin{split}
&\hskip-10pt|I_2(\lambda)|\\
\hskip-3pt\leq&\frac{1}{2\pi^{2}\rho_n^{2}}\sum_{v=1}^{n-1}\left|\frac{1}{n}\sum_{|u|<n}\!U_{n}(u,0,v)
~C\!\!\left(\frac{u}{\rho_n}\right)C\!\!\left(\frac{u\!+\!v}{\rho_n}\right)\!\right||K(b_n v)|\\
 +& \frac{1}{2\pi^{2}\rho_n^{2}}\sum_{v=1}^{n-1}
\left|\frac{1}{n}\sum_{|u|<n}\!U_{n}(u,0,v)~C\!\!\left(\frac{u}{\rho_n}\right)
C\!\!\left(\frac{u\!-\!v}{\rho_n}\right)\!\right||K(b_n v)|\\
+&\frac{1}{2\pi^{2}\rho_n^{2}}\sum_{v=1}^{n-1}\left|\frac{1}{n}\sum_{|u|<n}\!U_{n}(u,0,v)
~Q\!\!\left(\frac{v}{\rho_n},\frac{-u}{\rho_n},\frac{-u}{\rho_n}\right)\!\right||K(b_n v)|\\
\hskip-3pt\leq&\frac{|C(0)|}{\pi^{2}nb_n\rho_n}\left(\sum_{|u|<n}
\left|C\!\!\left(\frac{u}{\rho_n}\right)\right|\frac{1}{\rho_n}\right)\left(\sum_{v=1}^{n-1}|K(b_n v)|b_n\right)\\
&+
\frac{|g_1(0)g_2(0)|}{2\pi^{2}nb_n\rho_n}\left(\sum_{|u|<n}g_3\!\!
\left(\frac{u}{\rho_n}\right)\frac{1}{\rho_n}\right)\left(\sum_{v=1}^{n-1}|K(b_n
v)|b_n\right).
\end{split}
\end{equation*}
The last expression does not depend on $\lambda$. An argument as in
the case of $I_1$ will show that
$\lim_{n\rightarrow\infty}\sum_{v=1}^{n-1}|K(b_n v)|b_n=
\int_{0}^{\infty}K(x)dx$, under Conditions~2 and~3. The convergence
of the other sums have already been discussed in connection with the
term~$I_1$. Hence, $nb_nI_2(\lambda)\rightarrow 0$ as
$n\rightarrow\infty$ uniformly for all $\lambda$.

\medskip
Now we will consider $I_3(\lambda)$. Using (\ref{desofcov}), this
term can be written as
\begin{equation*}
\begin{split}
I_3(\lambda)=&\frac{1}{\pi^{2}\rho_n^{2}}\sum_{v_{1}=1}^{n-1}\sum_{v_{2}=1}^{n-1}
Cov[\widehat{\gamma}_n(v_{1}),\widehat{\gamma}_n(v_{2})]K(b_n
v_1)K(b_n v_2)\\&\times\cos\left(\frac{\lambda
v_1}{\rho_n}\right)\cos\left(\frac{\lambda
v_2}{\rho_n}\right)\\=&I_{31}(\lambda)+I_{32}(\lambda)+I_{33}(\lambda),
\end{split}
\end{equation*}
where
\begin{align*}
\begin{split}
I_{31}(\lambda)=&\frac{1}{\pi^{2}\rho_n^{2}}\sum_{v_1=1}^{n-1}\sum_{v_2=1}^{n-1}K(b_n
v_1)K(b_n v_2)\cos\!\left(\!\frac{\lambda
v_1}{\rho_n}\!\!\right)\cos\!\left(\!\frac{\lambda
v_2}{\rho_n}\!\!\right)\\
&\times\frac{1}{n}\sum_{u=-(n-1)}^{(n-1)}U_n(u,v_1,v_2)~C\!\left(\!\frac{u+v_1}{\rho_n}\!\right)
C\!\left(\!\frac{u-v_2}{\rho_n}\!\right)\!,
\end{split}
\end{align*}
\begin{align*}
\begin{split}
I_{32}(\lambda)=&\frac{1}{\pi^{2}\rho_n^{2}}\sum_{v_1=1}^{n-1}\sum_{v_2=1}^{n-1}K(b_n
v_1)K(b_n v_2)\cos\!\left(\!\frac{\lambda
v_1}{\rho_n}\!\!\right)\cos\!\left(\!\frac{\lambda v_2}{\rho_n}\!\!\right)\\
&\times\frac{1}{n}\sum_{u=-(n-1)}^{(n-1)}\!\!U_n(u,v_1,v_2)
~C\!\left(\!\frac{u}{\rho_n}\!\right)C\!\left(\!\frac{u+v_1-v_2}{\rho_n}\!\right)\!,
\end{split}
\end{align*}
and
\begin{align*}
\begin{split}
&\hskip-15pt I_{33}(\lambda)\\
=&\frac{1}{\pi^{2}\rho_n^{2}}\sum_{v_1=1}^{n-1}\sum_{v_2=1}^{n-1}\!K(b_n
v_1)K(b_n v_2)\cos\!\left(\!\frac{\lambda
v_1}{\rho_n}\!\right)\cos\!\left(\!\frac{\lambda
v_2}{\rho_n}\!\right)\!\\&\times
\frac{1}{n}\sum_{u=-(n-1)}^{(n-1)}U_n(u,v_1,v_2)~Q\!\left(\!\frac{v_1}{\rho_n},\frac{-u}{\rho_n},\frac{-u+v_2}{\rho_n}\!\!\right).
\end{split}
\end{align*}

From Conditions 2 and 5,
\begin{align*}
&\hskip-7pt |nb_nI_{33}(\lambda)|\\
\leq& M^{2}\rho_n
b_n\frac{1}{\rho_n^{3}}\!\sum_{v_1=1}^{n-1}\sum_{v_2=1}^{n-1}\sum_{u=-(n-1)}^{(n-1)}
\!\!g_1\!\!\left(\frac{v_1}{\rho_n}\!\right)g_2\!\!\left(\frac{u}{\rho_n}\!\right)g_3\!\!\left(\!\frac{u\!+\!v_2}{\rho_n}\!\right)\!.
\end{align*}
By using a similar argument as in the case of $I_1$, we have
\begin{align*}
&\hskip-5pt\lim_{n\rightarrow
\infty}\frac{1}{\rho_n^{3}}\!\sum_{v_1=1}^{n-1}\sum_{v_2=1}^{n-1}\sum_{u=-(n-1)}^{(n-1)}\!g_1\!\left(\frac{v_1}{\rho_n}\right)\!g_2\!\left(\frac{u}{\rho_n}\right)\!g_3\!\left(\!\frac{u\!+\!v_2}{\rho_n}\!\!\right)\\&=
\int_{0}^{\infty}
g_1(v_1)dv_1\int_{0}^{\infty}\left[\int_{-\infty}^{\infty}
g_2(u)g_3(u+v_2)du\right] dv_2.
\end{align*}
\medskip
Hence, $nb_nI_{33}(\lambda) \rightarrow 0$ as $n \rightarrow \infty$ uniformly for all $\lambda$.\\

Consider the term $I_{31}(\lambda)$, let
$j\!=\!v_1+v_2, ~ l\!=\!u-v_2 ,~ i\!=\!v_2$.
\begin{align*}
&\hskip-8pt I_{31}(\lambda)\\
=&\frac{1}{\pi^{2}
\rho_n^{2}}\sum_{i=1}^{n-1}\sum_{j=i+1}^{n-1+i}\!K(b_n i)K(b_n(
j\!-\! i))\cos\!\!\left(\frac{\lambda
i}{\rho_n}\!\right)\cos\!\!\left(\!\frac{\lambda(j\!-\!i)}{\rho_n}\!\right)\\
&\times\frac{1}{n}\sum_{l=-(n-1)-i}^{(n-1)-i}U_n(l+i,j-i,i)~C\!\!\left(\frac{j+l}{\rho_n}\right)C\!\!\left(\frac{l}{\rho_n}\right).
\end{align*}
From Condition 2, observe that
\begin{align*}
&\hskip-7pt nb_n|I_{31}(\lambda)|\\
&\leq
M\frac{b_n}{\rho_n^{2}}\sum_{i=1}^{n-1}\sum_{j=i+1}^{n-1+i}\sum_{l=-(n-1)-i}^{(n-1)-i}\left|K(b_ni)~C\!\!\left(\!\frac{j\!+\!l}{\rho_n}\!\right)C\!\!\left(\!\frac{l}{\rho_n}\!\right)\right|.
\end{align*}
Consider the simple function $S_n$, defined over
$(0,\infty)\times(0,\infty)\times(-\infty,\infty)$ by
\begin{align*}
&\hskip-7pt S_n(x,t,t^{'})\\
=&M\frac1{\pi^2}\sum_{i=1}^{n-1}\sum_{j=i+1}^{n-1+i}\sum_{l=-(n-1)-i}^{(n-1)-i}\left|K(b_n
i)~C\!\!\left(\!\frac{j\!+\!l}{\rho_n}\!\right)\!C\!\!\left(\!\frac{l}{\rho_n}\!\right)\right|\\
&\times\!1_{((i-1)b_n,i
b_n]}(x)1_{(\frac{j-1}{\rho_n},\frac{j}{\rho_n}]}(t)1_{(\frac{l-1}{\rho_n},\frac{l}{\rho_n}]}(t^{'}\!),
\end{align*}
so that
$$nb_n|I_{31}(\lambda)|\leq\int_{0}^{\infty}\int_{0}^{\infty}\int_{-\infty}^{\infty}
S_n(x,t,t^{'})dxdtdt^{'}.$$ %
Since $\lim_{n\rightarrow\infty}\rho_nb_n=0$, we have, for any fixed
$(x,t,t')\in(0,\infty)\times(0,\infty)\times(-\infty,\infty)$ and
for large enough $n$, the inequality $\rho_nb_n<x/t$, i.e.,
$t\rho_n<x/b_n$. Therefore, for large $n$, the unique integer $j$
for which $1_{(\frac{j-1}{\rho_n},\frac{j}{\rho_n}]}(t)$ is non-zero
is smaller than the unique integer $i$ for which $1_{((i-1)b_n,i
b_n]}(x)$ is non-zero. However, the ranges of summations in the
definition of $S_n(x,t,t^{'})$ do not permit the order $j\le i$.
Therefore, $\lim_{n\rightarrow\infty}S_n(x,t,t^{'})=0$. From
Condition 1 and 2, we have the dominance
$$|S_n(x,t,t^{'})|\leq M |K_1(x)h_0(t+t^{'})h_0(t^{'})| \in L^{1} $$
By applying the DCT, we have
$$\lim_{n\rightarrow\infty}\int_{0}^{\infty}\int_{0}^{\infty}\int_{-\infty}^{\infty}S_n(x,t,t^{'})dxdtdt^{'}=0.$$
Since the function $S_n(x,t,t^{'})$ does not depend on $\lambda$, we
have $nb_nI_{31}(\lambda)\rightarrow0$ uniformly for all $\lambda$.

\medskip

In view of the convergence of the terms $I_1$, $I_2(\lambda)$,
$I_{33}(\lambda)$ and $I_{31}(\lambda)$, we have
\begin{equation}\label{restvar}
\lim_{n\rightarrow\infty}nb_n\left[Var\left(\widehat{\phi}_{n}(\lambda)\right)-I_{32}(\lambda)\right]=0
\end{equation}
uniformly for all $\lambda$, and so we need to prove the convergence
of $nb_nI_{32}(\lambda)$ only.

Now consider $I_{32}(\lambda)$ and let $j=v_1-v_2 ,~ i=v_2,~ l=u$.
\begin{equation}\label{i32des}
\begin{split}
&\hskip-7pt
I_{32}(\lambda)\\=&\frac{1}{\pi^2\rho_n^2}\sum_{i=1}^{n-1}\sum_{j=-i+1}^{n-1-i}\!\!K(b_n
i)K(b_n (i\!+\! j))\cos\!\!\left(\frac{\lambda i}{\rho_n}\right)\cos\!\!\left(\!\frac{\lambda(i\!+\!j)}{\rho_n}\!\right)\\
&\times
\frac{1}{n}\sum_{l=-(n-1)}^{(n-1)}U_n(l,i+j,i)~C\!\!\left(\frac{l}{\rho_n}\right)C\!\!\left(\frac{l+j}{\rho_n}\right).
\end{split}
\end{equation}

For $\lambda=0$, it follows from (\ref{restvar}) and Lemma~1 below that
\begin{equation}
\lim_{n\rightarrow\infty}nb_nVar[\widehat{\phi}_{n}(0)]=2[\phi(0)]^{2}
\int_{-\infty}^{\infty}K^{2}(x)dx
\end{equation}

For $\lambda\ne 0$, we will further decompose
$I_{32}(\lambda)$ as follows. By applying the formula $2\cos(a)\cos(b)=\cos(a-b)+\cos(a+b)$ and $\cos(a+b)=\cos(a)\cos(b)-\sin(a)\sin(b)$, we have
$$I_{32}(\lambda)=I_{321}(\lambda)+I_{322}(\lambda)-I_{323}(\lambda),$$
where
\begin{equation}\label{a321}
\begin{split}
I_{321}(\lambda)=&\frac{1}{2\pi^2}\frac{1}{n\rho_n^2}\sum_{i=1}^{n-1}\sum_{j=-i+1}^{n-1-i}\!\!K(b_n
i)K(b_n( i\!+\! j))\cos\!\!\left(\frac{\lambda j}{\rho_n}\right)\\
&\times
\sum_{l=-(n-1)}^{(n-1)}U_n(l,i+j,i)~C\!\!\left(\frac{l}{\rho_n}\right)C\!\!\left(\frac{l+j}{\rho_n}\right),
\end{split}
\end{equation}
\begin{equation}\label{a322}
\begin{split}
&\hskip-7pt I_{322}(\lambda)\\
=&\frac{1}{2\pi^2}\frac{1}{n\rho_n^2}\sum_{i=1}^{n-1}\sum_{j=-i+1}^{n-1-i}\!\!K(b_n
i)K(b_n(i\!+\!j))\cos\!\!\left(\frac{\lambda
j}{\rho_n}\right)\cos\!\!\left(\frac{2\lambda  i}{\rho_n}\right)\\
&\times
\sum_{l=-(n-1)}^{(n-1)}U_n(l,i+j,i)~C\!\!\left(\frac{l}{\rho_n}\right)C\!\!\left(\frac{l+j}{\rho_n}\right),
\end{split}
\end{equation}
and
\begin{equation}\label{a323}
\begin{split}
&\hskip-7ptI_{323}(\lambda)\\
=&\frac{1}{2\pi^2}\frac{1}{n\rho_n^2}\sum_{i=1}^{n-1}\sum_{j=-i+1}^{n-1-i}K(b_n
i)K(b_n(i\!+\!j))\sin\!\!\left(\frac{\lambda
j}{\rho_n}\right)\sin\!\!\left(\frac{2\lambda  i}{\rho_n}\right)\\
&\times
\sum_{l=-(n-1)}^{(n-1)}U_n(l,i+j,i)~C\!\!\left(\frac{l}{\rho_n}\right)C\!\!\left(\frac{l+j}{\rho_n}\right).
\end{split}
\end{equation}

It follows from equation (\ref{restvar}), Lemma 2 and Lemma 3 below that
\begin{equation}
\begin{split}
\lim_{n\rightarrow\infty}nb_nVar[\widehat{\phi}_{n}(\lambda)]=
&\lim_{n\rightarrow\infty}nb_nI_{321}(\lambda)\\
=&[\phi(\lambda)]^{2} \int_{-\infty}^{\infty}K^{2}(x)dx,
\end{split}
\end{equation}
and the convergence is uniform over any closed interval that does
not include the frequency 0. This completes the proof.\\
\mbox{}\hfill$\Box$

\begin{Lem}
\begin{equation*}
\begin{split}
&\hskip-25pt\lim_{n\rightarrow\infty}
nb_nI_{32}(0)\\=&\frac{1}{\pi^2}\int_{0}^{\infty}\!\!K^{2}(x)\int_{-\infty}^{\infty}\left[\int_{-\infty}^{\infty}C(t+t^{'})C(t^{'})~dt^{'}\right]dt~dx.
\end{split}
\end{equation*}

\end{Lem}
\noindent {\bf Proof of Lemma 1.}\ \ Consider the simple function
$S_n(\cdot)$, defined over
$(0,\infty)\times(-\infty,\infty)\times(-\infty,\infty)$ by
\begin{align*}
&\hskip-7pt S_n(x,t,t^{'})\\
=&\frac{1}{\pi^2}\sum_{i=1}^{n-1}\sum_{j=-i+1}^{n-1-i}\sum_{l=-(n-1)}^{(n-1)}\!\!\!K(b_n i)K(b_n (i\!+\! j))
U_n(l,i\!+\!j,i)\\
&\times
\!C\!\!\left(\!\frac{l}{\rho_n}\!\right)\!C\!\!\left(\!\frac{l\!+\!j}{\rho_n}\!\right)\!1_{((i-1)b_n,i
b_n]}(x)1_{(\frac{j-1}{\rho_n},\frac{j}{\rho_n}]}(t)1_{(\frac{l-1}{\rho_n},\frac{l}{\rho_n}]}(t^{'}\!).
\end{align*}
Observe from (\ref{i32des}) that
$$nb_nI_{32}(0)=\int_{0}^{\infty}\int_{-\infty}^{\infty}\int_{-\infty}^{\infty}
S_n(x,t,t^{'})dxdtdt^{'}.$$ %
Define $i_n(x)$, $j_n(t)$ and $l_n(t^{'})$ as the smallest integers
greater than or equal to $x/b_n$, $\rho_nt$ and $\rho_nt'$,
respectively. Thus, $(x,t,t^{'}) \in (b_ni_{n-1}(x),b_n
i_n(x)]\times\left(\frac{j_{n-1}(t)}{\rho_n},\frac{j_n(t)}{\rho_n}\right]\times
\left(\frac{l_{n-1}(t')}{\rho_n},\frac{l_n(t')}{\rho_n}\right]$ and
$b_n i_n(x)\!\rightarrow\! x ,\frac{j_n(t)}{\rho_n}\!\rightarrow\!
t, \frac{l_n(t^{'})}{\rho_n}\! \rightarrow\! t^{'}$ as $n\rightarrow
\infty$. Since $nb_n\rightarrow\infty$ and $b_n\rho_n\rightarrow0$
as $n\rightarrow\infty$, we have, for any point $(x,t,t^{'})\in
(0,\infty)\times(-\infty,\infty)\times(-\infty,\infty)$ and large
enough $n$, the inequalities $-\frac
x{b_n\rho_n}<t<\frac{nb_n-x}{b_n\rho_n}{\displaystyle\phantom
\int}$\!\!\!\!\!, i.e., $-i_n(x)< j_n(t)< n-i_n(x)$. Thus, for
sufficiently large $n$, we have
\begin{align*}
S_n(x,t,t^{'})=&\frac{1}{\pi^2}K(b_n i_n(x))K\left(b_n
i_n(x)\!+\!b_n \rho_n \frac{j_n(t)}{\rho_n}
\right)\\&\times\!U_n(l_n(t^{'}\!),i_n(x)\!+\!j_n(t),i_n(x))\\&\times
C\left(\frac{l_n(t^{'})}{\rho_n}\right)C\left(\frac{l_n(t^{'})+j_n(t)}{\rho_n}\right).
\end{align*}
Also, for large $n$, we have $-(n\!-\!xb_n)+\frac
t{\rho_n}\!<\!\frac{t'}{\rho_n}\!<\!(n\!-\!xb_n)$, and so
$U_n(l_n(t^{'}\!),i_n(x)\!+\!j_n(t),i_n(x))$ is positive and it
converges to~1. Therefore, by virtue of Conditions~1 and~2, we have
$$\lim_{n\rightarrow\infty}S_n(x,t,t^{'})=\frac{1}{\pi^{2}}K^2(x)C(t^{'})C(t^{'}+t).$$
Again, from Conditions~1 and~2, we have the dominance
$$|S_n(x,t,t^{'})|\leq M |K_1(x)h_0(t+t^{'})h_0(t^{'})| \in L^{1}.
$$
By applying the DCT, we have
\begin{equation*}
\begin{split}
&\hskip-25pt \lim_{n\rightarrow\infty}nb_nI_{32}(0)\\
\!=\!&\lim_{n\rightarrow\infty}\int_{0}^{\infty}\int_{-\infty}^{\infty}\int_{-\infty}^{\infty}
S_n(x,t,t^{'})dxdtdt'\\\!=\!&\frac{1}{\pi^2}\!\int_{0}^{\infty}\!\!\!K^{2}(x)\!\int_{-\infty}^{\infty}\!\!\left[\int_{-\infty}^{\infty}\!\!\!C(t+t^{'}\!)C(t^{'}\!)dt^{'}\!\right]\!dtdx.
\end{split}
\end{equation*}
\hfill$\Box$

\begin{Lem}
The function $I_{321}(\cdot)$ converges as follows:
\begin{equation*}
\begin{split}
&\hskip-10pt\lim_{n\rightarrow \infty}nb_nI_{321}(\lambda)\\
=&\frac{1}{2\pi^{2}}\int_{0}^{\infty}\!\!\!
K^{2}(x)\int_{-\infty}^{\infty}\!\! \left[\cos(\lambda
t)\!\left[\int_{-\infty}^{\infty}\!\!C(t+t^{'}\!)C(t^{'}\!)dt^{'}\right]\!dt\right]\!dx.
\end{split}
\end{equation*}
The convergence is uniform on $[\lambda_l,\lambda_u]$ for arbitrary
$\lambda_l$ and $\lambda_u$ such that $\lambda_l<\lambda_u$ and
$\lambda_l\lambda_u>0$.
\end{Lem}
\noindent{\bf Proof of Lemma 2.}\ \ Consider the simple function
$S_n(\cdot)$, defined over
$[\lambda_l,\lambda_u]\times(0,\infty)\times(-\infty,\infty)\times(-\infty,\infty)$
by
\begin{align*}
S_n(\lambda,x,t,t^{'})
=&\frac{1}{2\pi^2}\!\sum_{i=1}^{n-1}\sum_{j=-i+1}^{n-1-i}\!\sum_{l=-(n-1)}^{(n-1)}\!\!\!K(b_n
i)K(b_n (i\!+\! j))\\
&\times\cos\!\!\left(\!\frac{\lambda
j}{\rho_n}\!\right)U_n(l,i\!+\!j,i)~C\!\!\left(\frac{l}{\rho_n}\!\right)C\!\!\left(\frac{l\!+\!j}{\rho_n}\right)\\
&\times 1_{((i-1)b_n,i
b_n]}(x)1_{(\frac{j-1}{\rho_n},\frac{j}{\rho_n}]}(t)1_{(\frac{l-1}{\rho_n},\frac{l}{\rho_n}]}(t^{'}\!),
\end{align*}
so that, from (\ref{a321}),
$$nb_nI_{321}(\lambda)=\int_{0}^{\infty}\int_{-\infty}^{\infty}\int_{-\infty}^{\infty}
S_n(\lambda,x,t,t^{'})dxdtdt^{'}.$$ %
A similar argument as in the proof of Lemma 1 will show that for
$(x,t,t^{'})\in
(0,\infty)\times(-\infty,\infty)\times(-\infty,\infty)$ and
sufficiently large $n$,
\begin{align*}
S_n(\lambda,x,t,t^{'}) =&\frac{1}{2\pi^2}K(b_n
i_n(x))~K\!\!\left(b_n i_n(x)+b_n\rho_n
\frac{j_n(t)}{\rho_n}\right)\\
&\times\cos\!\!\left(\!\frac{\lambda
j_n(t)}{\rho_n}\!\!\right) U_n(l_n(t^{'}),i_n(x)\!+\!j_n(t),i_n(x))\\
&\times
C\!\!\left(\frac{l_n(t^{'})}{\rho_n}\right)C\!\!\left(\frac{l_n(t^{'})\!+\!j_n(t)}{\rho_n}\right),
\end{align*}
where $i_n(x)$, $j_n(t)$ and $l_n(t^{'})$ are the smallest integers
greater than or equal to $x/b_n$, $\rho_nt$ and $\rho_nt'$,
respectively, and that the function $S_n(\lambda,x,t,t^{'})$
converges to the function $S(\cdot)$, defined over
$(0,\infty)\times(-\infty,\infty)\times(-\infty,\infty)$ by
$$S(\lambda,x,t,t^{'})=\frac{1}{2\pi^2}K^{2}(x)\cos(\lambda t)C(t^{'})C(t^{'}+t).$$

Observe also that
$\int_{0}^{\infty}\int_{-\infty}^{\infty}\int_{-\infty}^{\infty}S(\lambda,x,t,t^{'})dxdtdt^{'}$
is a continuous function in $\lambda$. As in the proof of Theorem~1,
we prove the convergence of
$\int_{0}^{\infty}\int_{-\infty}^{\infty}\int_{-\infty}^{\infty}S_n(\lambda,x,t,t^{'})dxdtdt^{'}$%
uniformly on $[\lambda_l,\lambda_u]$, by showing that for any
sequence $\lambda_n \rightarrow \lambda$,
\begin{equation*}
\begin{split}
&\hskip-40pt \lim_{n\rightarrow\infty}
\int_{0}^{\infty}\int_{-\infty}^{\infty}\int_{-\infty}^{\infty}S_n(\lambda_n,x,t,t^{'})dxdtdt^{'}\\
&\mbox{}\qquad=\int_{0}^{\infty}\int_{-\infty}^{\infty}\int_{-\infty}^{\infty}S(\lambda,x,t,t^{'})dxdtdt^{'}
\end{split}
\end{equation*}
for $\lambda_n,\lambda \in [\lambda_l,\lambda_u] $. The latter
convergence follows, through Condition~1 and~2 and the DCT, from the dominance
$$|S_n(\lambda,x,t,t^{'})|\leq M |K_1(x)h_0(t+t^{'})h_0(t^{'})| \in L^{1}, $$
and the convergence of the integrand, which holds because of the
continuity of the kernel, the cosine and the covariance function.

Hence, $nb_nI_{321}(\cdot)$ converges as stated uniformly on
$[\lambda_l,\lambda_u]$.\\
\mbox{}\hfill$\Box$

\begin{Lem}
The functions $nb_nI_{322}(\cdot)$ and $nb_nI_{323}(\cdot)$ converge to $0$
uniformly on $[\lambda_l,\lambda_u]$ for arbitrary $\lambda_l$ and
$\lambda_u$ such that $\lambda_l<\lambda_u$ and
$\lambda_l\lambda_u>0$.
\end{Lem}
\noindent {\bf Proof of Lemma 3.}\ \ Consider the simple function
$S_n(\cdot)$, defined over
$[\lambda_l,\lambda_u]\times(0,\infty)\times(-\infty,\infty)\times(-\infty,\infty)$
by
\begin{align*}
&\hskip-10pt
S_n(\lambda,x,t,t^{'})\\=&\frac{1}{2\pi^2}\sum_{i=1}^{n-1}\sum_{j=-i+1}^{n-1+i}\sum_{l=-(n-1)}^{(n-1)}K(b_n
i)K(b_n( i+j))\cos\!\left(\frac{\lambda j}{\rho_n}\right)\\
&\times\cos\!\left(\frac{2\lambda b_n i}{\rho_n b_n}\right)
U_n(l,i\!+\!j,i)~C\!\!\left(\frac{l}{\rho_n}\right)C\!\!\left(\frac{l+j}{\rho_n}\right)\\&\times1_{((i-1)b_n,i
b_n]}(x)\times1_{(\frac{j-1}{\rho_n},\frac{j}{\rho_n}]}(t)1_{(\frac{l-1}{\rho_n},\frac{l}{\rho_n}]}(t^{'}),
\end{align*}
so that, from (\ref{a322}),
\begin{equation}
nb_nI_{322}(\lambda)=\int_{0}^{\infty}\int_{-\infty}^{\infty}\int_{-\infty}^{\infty}
S_n(\lambda,x,t,t^{'})dxdtdt^{'}. \label{i322}
\end{equation} %
A similar argument as in the proof of Lemma 1 will show that for
$(x,t,t^{'})\in
(0,\infty)\times(-\infty,\infty)\times(-\infty,\infty)$ and
sufficiently large $n$,
\begin{align*}
&\hskip-10pt S_n(\lambda,x,t,t^{'})\\=&\frac{1}{2\pi^2}K(b_n
i_n(x))K\left(b_ni_n(x)+b_n\rho_n\frac{j_n(t)}{\rho_n}\right)\cos\!\left(\frac{\lambda
j_n(t)}{\rho_n}\right)\\
&\times\cos\!\left(\frac{2\lambda i_n(x)b_n}{b_n\rho_n}\right)U_n(l_n(t^{'}),j_n(t)+i_n(x),i_n(x))\\
&\times
C\!\!\left(\frac{l_n(t^{'})}{\rho_n}\right)C\!\!\left(\frac{l_n(t^{'})+j_n(t)}{\rho_n}\right),
\end{align*}
where $i_n(x)$, $j_n(t)$ and $l_n(t^{'})$ are the smallest integers
greater than or equal to $x/b_n$, $\rho_nt$ and $\rho_nt'$,
respectively.

For obtaining the uniform convergence of $I_{322}(\cdot)$, consider
\begin{align}
&\hskip-10pt\sup_{\lambda \in
[\lambda_l,\lambda_u]}\left|\int_{0}^{\infty}\!\!\int_{-\infty}^{\infty}\!\!\int_{-\infty}^{\infty}\!\!S_n(\lambda,x,t,t^{'})dxdtdt^{'}\right|\notag\\
\!\le\!&\sup_{\lambda \in
[\lambda_l,\lambda_u]}\int_{0}^{\infty}\!\!\int_{-\infty}^{\infty}\!\!\int_{-\infty}^{\infty}\!\!|S_n(\lambda,x,t,t^{'}\!)-g_n(\lambda,x,t,t^{'}\!)|dxdtdt^{'}\notag\\
~+\!&\sup_{\lambda \in
[\lambda_l,\lambda_u]}\left|\int_{0}^{\infty}\int_{-\infty}^{\infty}\int_{-\infty}^{\infty}g_n(\lambda,x,t,t^{'})dxdtdt^{'}\right|,
\label{2step}
\end{align}
where the function $g_n(\cdot)$ is defined over
$[\lambda_l,\lambda_u]\times(0,\infty)\times(-\infty,\infty)\times(-\infty,\infty)$ by
$$g_n(\lambda,x,t,t^{'})
=\frac{1}{2\pi^2}\cos(\lambda t)\cos\left(\frac{2\lambda
x}{b_n\rho_n}\right) K^2(x)C(t^{'})C(t+t^{'}).$$ %
We shall prove the convergence of $nb_nI_{322}(\cdot)$ given in
(\ref{i322}) by proving the convergence of the two integrals on the
right hand side of (\ref{2step}).

In order to prove the first convergence, we follow the route taken
in Theorem~1, i.e., we show that for any sequence
$\lambda_n\!\rightarrow\!\lambda$
$$\lim_{n\rightarrow\infty}\int_{0}^{\infty}\!\!\!\int_{-\infty}^{\infty}\!\!\int_{-\infty}^{\infty}\!\!|S_n(\lambda_n,x,t,t^{'})-g_n(\lambda_n,x,t,t^{'})|dxdtdt^{'}=0,$$
for $\lambda,\lambda_n\in[\lambda_l,\lambda_u]$. The above integral
can be written as
\begin{align}
&\hskip-10pt
\int_{0}^{\infty}\int_{-\infty}^{\infty}\int_{-\infty}^{\infty}|S_n(\lambda_n,x,t,t^{'})-g_n(\lambda_n,x,t,t^{'})|dxdtdt'\notag\\
\le&\int_{0}^{\infty}\int_{-\infty}^{\infty}\int_{-\infty}^{\infty}|S_n(\lambda_n,x,t,t^{'})-G_n(\lambda_n,x,t,t^{'})|dxdtdt'\notag\\
&+\int_{0}^{\infty}\int_{-\infty}^{\infty}\int_{-\infty}^{\infty}|G_n(\lambda_n,x,t,t^{'})-g_n(\lambda_n,x,t,t^{'})|dxdtdt',
\label{2part}
\end{align}
where the function $G_n(\cdot)$ is defined over
$[\lambda_l,\lambda_u]\times(0,\infty)\times(-\infty,\infty)\times(-\infty,\infty)$ by
\begin{align*}
&\hskip-15pt G_n(\lambda,x,t,t^{'})&\\
=&\frac{1}{2\pi^2} \cos(\lambda t)\cos\!\left(\!\frac{2\lambda
i_n(x) b_n}{b_n\rho_n}\right) K^2(x)C(t^{'})C(t+t^{'})).
\end{align*}
Now observe that
\begin{align*}
&\hskip-45pt|S_n(\lambda_n,x,t,t^{'})-G_n(\lambda_n,x,t,t^{'})|\\\le&
M \left|\cos\left(\frac{2\lambda_n i_n(x)
b_n}{b_n\rho_n}\right)\alpha_n(\lambda_n,x,t,t^{'})\right|,
\end{align*}
where
\begin{align*}
&\hskip-8pt\alpha_n(\lambda_n,x,t,t^{'})\\
=&K(b_n i_n(x))\cos\!\left(\frac{\lambda_n
j_n(t)}{\rho_n}\right)U_n(l_n(t^{'}),j_n(t)\!+\!i_n(x),i_n(x))\\
&\times C\!\!\left(\frac{l_n(t^{'})}{\rho_n}\right)C\!\!\left(\frac{l_n(t^{'})+j_n(t)}{\rho_n}\right)\\
&-\cos(\lambda_n t)K(x)C(t^{'})C(t+t^{'}).
\end{align*}
Since $\alpha_n(\lambda_n,x,t,t^{'})\rightarrow 0$ as $n\rightarrow
\infty$, we have
$$\lim_{n\rightarrow\infty}|S_n(\lambda_n,x,t,t^{'})-G_n(\lambda_n,x,t,t^{'})|= 0.$$
Since from Condition 1 and 2, we have the dominance
\begin{align*}
&\hskip-40pt |S_n(\lambda_n,x,t,t^{'})-G_n(\lambda_n,x,t,t^{'})|\\&\le 2M|K_1(x)h_0(t')h_0(t+t')|\in L^{1}.
\end{align*}
By applying DCT, we have
$$\lim_{n\rightarrow\infty}\int_{0}^{\infty}\!\!\!\!\int_{-\infty}^{\infty}\!\!\!
\int_{-\infty}^{\infty}\!\!|S_n(\lambda_n,x,t,t^{'})-G_n(\lambda_n,x,t,t^{'})|dxdtdt'\!=\!0.$$

Turning to the second term on the right hand side of (\ref{2part}),
observe that for any fixed $x$, $|x- i_n(x)b_n|\le b_n$. By applying
the Mean Value Theorem to the cosine function in the interval
$[\frac{\lambda_n x}{\rho_n b_n}, \frac{\lambda_n i_n(x)b_n}{\rho_n
b_n}]$, we have
\begin{align*}
&\hskip-50pt \cos\left(\frac{\lambda_n i_n(x)b_n}{\rho_n b_n}\right)-\cos\left(\frac{\lambda_n x}{\rho_n b_n}\right)\\
=&-\sin(\theta)\left|\frac{\lambda_ni_n(x)b_n}{\rho_nb_n}-\frac{\lambda_nx}{\rho_nb_n}\right|,
\end{align*}
for some $\theta \in \left[\frac{\lambda_n x}{\rho_n b_n}, \frac{\lambda_n i_n(x)b_n}{\rho_n b_n}\right]$.
Therefore
$$\left|\cos\left(\!\!\!\frac{\lambda_n i_n(x)b_n}{\rho_n b_n}\!\right)-\cos\left(\!\frac{\lambda_n x}{\rho_n b_n}\!\right)\right|\le \frac{\lambda_n}{\rho_n}.$$
Thus,
\begin{align*}
&\hskip-15pt|G_n(\lambda_n,x,t,t^{'})-g_n(\lambda_n,x,t,t^{'})|\\
\leq& M^{2}C^{2}(0)\left|\cos\left(\frac{\lambda_n i_n(x)b_n}{\rho_n b_n}\right)
-\cos\left(\frac{\lambda_n x}{\rho_n b_n}\right)\right|\\
\le& M^{2}C^{2}(0)  \frac{\lambda_n}{\rho_n}.
\end{align*}
So
$$\lim_{n\rightarrow\infty}|G_n(\lambda_n,x,t,t^{'})-g_n(\lambda_n,x,t,t^{'})|\rightarrow 0.$$
From Condition 1 and 2, we have the dominance
\begin{align*}
&\hskip-50pt|G_n(\lambda_n,x,t,t^{'})-g_n(\lambda_n,x,t,t^{'})|\\
\leq&2M|K_1(x)h_0(t+t^{'})h_0(t^{'})|\in L^{1},
\end{align*}
which leads us, through another use of the DCT, the convergence of
the second integral of (\ref{2part}). This establishes that the
first term on the right hand side of (\ref{2step}) converges to~0. We only have to deal with the second term.

Let
$$s_n(\lambda)=\int_{0}^{\infty}\int_{-\infty}^{\infty}\int_{-\infty}^{\infty}g_n(\lambda,x,t,t^{'})dxdtdt^{'}.$$
In order to establish the uniform convergence of $s_n(\cdot)$ over
$[\lambda_l,\lambda_u]$, it is enough to show that
$s_n(\lambda_n)\rightarrow 0$ for any sequence
$\lambda_n\rightarrow\lambda$, where
$\lambda,\lambda_n\in[\lambda_l,\lambda_u]$. By using the
Reimann-Lebesgue lemma, we have $s_n(\lambda_n)\rightarrow 0$. Thus,
the second term on the right hand side of (\ref{2step}) also
converges to~0. Hence, $nb_nI_{322}(\lambda)$ converges to~0
uniformly on $[\lambda_l,\lambda_u]$ as $n\rightarrow\infty$.

Convergence of $nb_nI_{323}(\lambda)$ to 0 can be established in a
similar manner.$\hfill\Box$

\bigskip \noindent{\bf Proof of Theorem 3.}\ \ Condition~1 ensures
absolute summability of the covariance sequence $\{C(k/\rho_n),~k\in
\mathbb{Z}\}$ of the regularly sampled process
$\left\{X\left(\frac{t}{\rho_{n}}\right) ,\ t \in \mathbb{Z}\right\}
$, for fixed $\rho_{n}$. The corresponding spectral density
$\xi_n(\cdot)$ is defined as
$$ \xi_n(\lambda)=\frac{1}{2\pi\rho_n}\sum_{j=-\infty}^{\infty}C\left(\frac{j}{\rho_n}\right)e^{-\frac{ij\lambda}{\rho_n}},~~\lambda\in(-\infty,\infty).$$
The function $\xi_n(\cdot)$ is periodic with period $2\pi\rho_n$ and
is related to the function $\phi(\cdot)$ as follows:
$$\xi_n(\lambda)=\sum_{l=-\infty}^{\infty}\phi(\lambda+2\pi l \rho_n), \qquad\lambda\in(-\pi\rho_n,\pi\rho_n].$$
In particular, for $\lambda\in [-\pi\rho_n,\pi\rho_n]$,
$$ \phi(\lambda)=\xi_n(\lambda)-\sum_{l=1}^{\infty}\phi(\lambda+2\pi l \rho_n)-\sum_{l=-\infty}^{-1}\phi(\lambda+2\pi l \rho_n).$$

For sufficiently large $n$, $\pi\rho_n$ lies outside any finite
interval $[\lambda_l,\lambda_u]$, and the bias of the estimator
$\widehat{\phi}_n(\lambda)$ given by (\ref{unifestn}) on
$[\lambda_l,\lambda_u]$ can be decomposed as follows.
\begin{equation}
\begin{split}
&\hskip-8pt E[\widehat{\phi}_n(\lambda)]-\phi(\lambda)\\=&
\frac{1}{2\pi\rho_n}\sum_{|v|<n}\left(1-\frac{|v|}{n}\right)C\left(\frac{v}{\rho_n}\right)K(b_n
v)e^{-\frac{iv\lambda}{\rho_n}}-\phi(\lambda)\\=&\frac{1}{2\pi\rho_n}\sum_{|v|<n}\left(1-\frac{|v|}{n}\right)C\left(\frac{v}{\rho_n}\right)K(b_n
v)e^{-\frac{iv\lambda}{\rho_n}}-\xi_n(\lambda)\\&+\sum_{l=1}^{\infty}\phi(\lambda+2\pi
l \rho_n)+\sum_{l=-\infty}^{-1}\phi(\lambda+2\pi l
\rho_n)\\=&B_1(\lambda)+B_2(\lambda)+B_3(\lambda)+B_4(\lambda)+B_5(\lambda),
\end{split}
\end{equation}
where
\begin{IEEEeqnarray*}{rCl}
B_1(\lambda)&=&-\frac{1}{2\pi\rho_n}\sum_{|v|<n}(1-K(b_n
v))C\left(\frac{v}{\rho_n}\right)e^{-\frac{i v
\lambda}{\rho_n}},\\
B_2(\lambda)&=&-\frac{1}{2\pi\rho_n}\sum_{|v|<n}\frac{|v|}{n}C\left(\frac{v}{\rho_n}\right)K(b_n
v)e^{-\frac{iv\lambda}{\rho_n}},\\
B_3(\lambda)&=&-\frac{1}{2\pi\rho_n}\sum_{|v|\ge
n}C\left(\frac{v}{\rho_n}\right)e^{-\frac{iv\lambda}{\rho_n}},\\
B_4(\lambda)&=&\sum_{l=1}^{\infty}\phi(\lambda+2\pi
l \rho_n),\\
B_5(\lambda)&=&\sum_{l=-\infty}^{-1}\phi(\lambda+2\pi l
\rho_n).
\end{IEEEeqnarray*}
We will consider each $B_i(\lambda)$, $i=1,\ldots,5$, separately.
\begin{align*}
&\hskip-25pt \left(\frac{1}{\rho_nb_n}\!\right)^{q}B_1(\lambda)\\
=&-\frac{1}{2\pi}\!\!\sum_{|v|<n}\left(\frac{1-K(b_n v)}{b_n^q
|v|^q}\!\right)\frac{|v|^q}{\rho_n^q}\!C\left(\frac{v}{\rho_n}\right)
e^{-\frac{iv\lambda}{\rho_n}}\frac{1}{\rho_n}.
\end{align*}
Consider the simple function $S_n$ defined over
$[\lambda_l,\lambda_u]\times(-\infty,\infty)$ as
\begin{align*}
&\hskip-8pt S_n(\lambda,t)\\
=&-\frac{1}{2\pi}\sum_{|v|<n}\!\!\left(\frac{1-K(b_n v)}{b_n^q
|v|^q}\right)\!\frac{|v|^q}{\rho_n^q}\,C\!\!\left(\frac{v}{\rho_n}\right)
e^{-\frac{iv\lambda}{\rho_n}}1_{(\frac{v-1}{\rho_n},\frac{v}{\rho_n}]}(t).
\end{align*}
Observe that $\left(\frac{1}{\rho_n
b_n}\right)^{q}B_1(\lambda)=\int_{-\infty}^{\infty}S_n(\lambda,t)dt$.

For any $t\in(-\infty,\infty)$, we define $v_n(t)$ as the smallest
integer greater than or equal to $t\rho_n$. It follows that
$\frac{v_n(t)}{\rho_n} \rightarrow t ~\mbox{as}~ n\rightarrow
\infty$, and for sufficiently large $n$ and any
$t\in(-\infty,\infty)$, we can write
$$S_n(\lambda,t)\!=\!-\frac{1}{2\pi}\left(\!\frac{1\!-\!K(b_n v_n(t))}{b_n^q
|v_n(t)|^q}\!\right)\!\frac{|v_n(t)|^q}{\rho_n^q}\,C\!\!\left(\!\frac{v_n(t)}{\rho_n}\!\right)
e^{-\frac{iv_n(t)\lambda}{\rho_n}}.
$$
From Conditions~4 and~2A, we have
$$\lim_{n\rightarrow \infty}S_n(t)=-\frac{1}{2\pi}k_{q}|t|^{q}C(t)e^{-it\lambda}$$
Also, Condition 1A implies
$$|S_n(t)|\leq M_1 h_q(t)~ \mbox{where} ~M_1=\sup_{x}\frac{1-K(x)}{x^{q}}.$$
By applying the DCT, we have
\begin{align}
\lim_{n\rightarrow \infty}\left(\frac{1}{\rho_n b_n}\right)^{q}B_1(\lambda)
\!=\!&\lim_{n\rightarrow \infty}\int_{-\infty}^{\infty}S_n(\lambda,t)dt\notag\\
\!=\!&-k_{q}\frac{1}{2\pi}\int_{-\infty}^{\infty}|t|^{q}C(t)e^{-it\lambda}
dt.\label{2arelax}
\end{align}
Thus, $B_1(\lambda)$ is $O\left((\rho_nb_n)^q\right)$. The fact that
this convergence is uniform over the interval
$[\lambda_l,\lambda_u]$ can be established by choosing any sequence
$\lambda_n$ in this interval that converges to $\lambda$, and
showing that $\int_{-\infty}^{\infty}S_n(\lambda_n,t)dt$ converges
to the right hand side of (\ref{2arelax}).

\medskip The term $B_2(\lambda)$ can be written as
$$\frac{n}{\rho_n}B_2(\lambda)=\int_{-\infty}^{\infty}S_n(\lambda,t)dt,$$
where $S_n(\cdot)$ is defined over
$[\lambda_l,\lambda_u]\times(-\infty,\infty)$ as
$$S_n(\lambda,t)=-\frac{1}{2\pi}\sum_{|v|<n}\frac{|v|}{\rho_n}C\left(\frac{v}{\rho_n}\right)K(b_n
v)e^{-\frac{iv\lambda}{\rho_n}}1_{(\frac{v-1}{\rho_n},\frac{v}{\rho_n}]}(t).$$
As in the case of $B_1(\lambda)$, it can be shown that
$$\lim_{n\rightarrow\infty}S_n(\lambda,t)=-\frac{1}{2\pi}|t|C(t)e^{-it\lambda}.$$
From Condition 1A, it follows that $|S_n(\lambda,t)|\le h_q(t) $.
Now again by applying the DCT, we have %
\begin{align*}
\lim_{n\rightarrow\infty}\frac{n}{\rho_n}B_2(\lambda)=&\lim_{n\rightarrow\infty}\int_{-\infty}^{\infty}S_n(\lambda,t)dt\\
=&-\frac{1}{2\pi}\int_{-\infty}^{\infty}|t|C(t)e^{-it\lambda}dt.
\end{align*}
Thus, $B_2(\lambda)$ is $O(\frac{\rho_n}{n})$. The uniform
convergence can be argued similarly as in the case of
$B_1(\lambda)$.

\medskip
The term $B_3(\lambda)$ satisfies
\begin{align*}
\hskip-15pt|B_3(\lambda)|\le&\frac{1}{2\pi\rho_n}\sum_{|v|\ge
n}\left|C\left(\frac{v}{\rho_n}\right)\right|
\\=&\left(\frac{\rho_n}{n}\right)^{q}\frac{1}{2\pi}\sum_{|v|\ge
n}\left(\frac{n}{\rho_n}\right)^q\left|C\left(\frac{v}{\rho_n}\right)\frac{1}{\rho_n}\right|\\
\le& \left(\frac{\rho_n}{n}\right)^{q}\frac{1}{\pi}\left[\sum_{v\ge
1}\left(\frac{v}{\rho_n}\right)^q\left|C\left(\frac{v}{\rho_n}\right)\frac{1}{\rho_n}\right|\right].
\end{align*}
Observe that for each fixed $n$, we have from Condition 1A,
\begin{equation*}
\begin{split}
&\hskip-45pt \lim_{m \rightarrow \infty}\sum_{v=1}^{m}\left(\frac{v}{\rho_n}\right)^q\left|C\left(\frac{v}{\rho_n}\right)\right|\frac{1}{\rho_n}\\&\le\lim_{m \rightarrow \infty}\sum_{v=1}^{m}h_q\left(\frac{v}{\rho_n}\right)\frac{1}{\rho_n}\\&\le \lim_{m \rightarrow \infty}\int_{0}^{\frac{m}{\rho_n}}h_q(t)dt\le\int_{0}^{\infty}h_q(t)dt.
\end{split}
\end{equation*}
So
$$
\limsup_{n\rightarrow
\infty}\sum_{v=1}^{\infty}\left(\frac{v}{\rho_n}\right)^q\left|C\left(\frac{v}{\rho_n}\right)\right|\frac{1}{\rho_n}\le\int_{0}^{\infty}h_q(t)dt.$$
Hence $|B_3(\lambda)|$ is bounded by an $O((\frac{\rho_n}{n})^q)$
term, which converges to zero faster than $B_2(\lambda)$.

\medskip
As for the term $B_4(\lambda)$, we have from Condition 1B and DCT
\begin{align*}
&\hskip-10pt\lim_{n\rightarrow\infty}(\rho_n)^{p}B_4(\lambda)\\
=&\lim_{n\rightarrow\infty}\sum_{l=1}^{\infty}\frac{(\rho_n)^p}{|\lambda\!+\!2\pi
l\rho_n|^p} |\lambda\!+\!2\pi l\rho_n|^p \phi(\lambda\!+\!2\pi l\rho_n)\\
=&\frac 1{(2\pi)^p}
\sum_{l=1}^{\infty}\lim_{n\rightarrow\infty}\frac{1}{|\frac{\lambda}{2\pi\rho_n}\!+\!l|^p}
\!\left(\lim_{n\rightarrow\infty}|\lambda\!+\!2\pi l\rho_n|^p
\phi(\lambda\!+\!2\pi l\rho_n)\!\right)\\
=&\frac{A}{(2\pi)^{p}}\sum_{l=1}^{\infty}\frac{1}{|l|^p}.
\end{align*}
Hence $B_4(\lambda)=O\left(\frac{1}{(\rho_n)^{p}}\right)$.

\medskip
Similarly it can be proved that
$B_5(\lambda)=O\left(\frac{1}{(\rho_n)^{p}}\right).$

\medskip
The theorem is proved by combining the five terms. \hfill$\Box$

\bigskip \noindent{\bf Proof of Theorem 4.}\ \ It follows from
Theorems~2 and~3 that the MSE of the estimator
$\widehat\phi_n(\cdot)$ can be written as
\begin{eqnarray}
MSE[\widehat{\phi}_n(\lambda)]&=&[E\{\widehat\phi_n(\lambda)-\phi(\lambda)\}]^2+
Var[\widehat{\phi}_n(\lambda)]\nonumber\\
&=&O\left((\rho_nb_n)^{2q}\right)+O\left({\rho_n^2\over
n^2}\right)+O\left({1\over\rho_n^{2p}}\right)\nonumber\\
&&\qquad+O\left({1\over
nb_n}\right). %
\label{mserate}
\end{eqnarray}
Let us first fix $n$ and $\rho_n$ and minimize the MSE with respect
to $b_n$. The squared bias is an increasing functions of $b_n$,
while the variance is a decreasing function of $b_n$. Therefore, the
maximum possible value is minimized (i.e., the fastest rate of
convergence is achieved) when $(\rho_nb_n)^{2q}\propto(nb_n)^{-1}$,
i.e., when
\begin{equation}
b_n\propto\left(n\rho_n^{2q}\right)^{-{1\over 2q+1}}. \label{interm}
\end{equation}
By substituting this value in the expression for the MSE, and making
use of the fact that ${2q\over 2q+1}<2$ and $\rho_n/n<1$, we have
$$MSE[\widehat{\phi}_n(\lambda)]=O\left(\left({\rho_n\over n}\right)^{2q\over 2q+1}\right)+O\left({1\over\rho_n^{2p}}\right).$$
The first term on the right hand side is an increasing function of
$\rho_n$, while the second term is a decreasing function of
$\rho_n$. Therefore, the maximum of the two terms is minimized when
$\left({\rho_n\over n}\right)^{2q\over 2q+1}\propto\rho_n^{-2p}$,
i.e., when $\rho_n$ is chosen as in (\ref{rhon}). The optimal rate
for $b_n$, as given in (\ref{bn}), is obtained by substituting the
expression for $\rho_n$ in (\ref{interm}). Further substitution of
these two optimal rates in (\ref{mserate}) gives
(\ref{msereg}).\hfill$\Box$

\section*{Acknowledgements}
The authors gratefully acknowledge suggestions and technical help
from Professors B.V. Rao and Arup Bose of the Indian Statistical
Institute. Suggestions from two referees have been useful in
improving the clarity of the presentation.

\begin{IEEEbiographynophoto}{Radhendushka Srivastava}
received the B.Sc. degree and the M.Sc. degree in statistics from
the University of Lucknow in the years 2003 and 2005, respectively.
He is currently a senior research fellow in the Indian Statistical
Institute, and is working towards the Ph.D. degree in statistics.
His research interests include time series analysis, stochastic
processes and applications of statistical methods to signal
processing.
\end{IEEEbiographynophoto}

\begin{IEEEbiographynophoto}{Debasis Sengupta}
received B.Tech. degree in electronics and electrical communications
engineering from the Indian Institute of Technology, Kharagpur, in
1984 and master's degrees in electrical engineering and in
statistics from the University of Rhode Island, Kingston in the year
1986 and the University of California, Santa Barbara, in the year
1988, respectively. He received Ph.D. degrees in electrical and
computer engineering and in Statistics from the University of
California, Santa Barbara, in the years 1989 and 1990, respectively.
Subsequently he has been with the Applied Statistics Unit of the
Indian Statistical Institute, Kolkata, where he has been a professor
since 1997. His research interests include regression, multivariate
analysis, time series analysis and statistical signal processing.
\end{IEEEbiographynophoto}

\end{document}